\documentclass[12pt,leqno,fleqn]{article}
\usepackage{epsf,amsmath,amssymb,amscd,hyperref,srcltx}
\usepackage{amsfonts}
\usepackage{psfrag}
\newcommand{\dps}{\displaystyle}
\newcommand{\GL}{\text{\rm GL}}
\newcommand{\tussenkop}[1]{\section*{\hspace*{\fill}\normalsize\boldmath #1\hspace*{\fill}}}
\newcommand{\GG}{{\cal G}}
\newcommand{\OO}{{\cal O}}
\newcommand{\RR}{{\cal R}}
\newcommand{\WW}{{\cal W}}
\newcommand{\inn}{^{\text{\rm in}}}
\newcommand{\uit}{^{\text{\rm out}}}
\newcounter{bewering}
\newcommand{\dy}[2]{%
\refstepcounter{equation}%
\LABEL{#1}%
\begin{list}{}{
\topsep 3mm
\leftmargin 18mm
\rightmargin 0cm
\itemsep 0mm
\listparindent 0mm
\parsep 0mm
\itemsep 0mm
\labelsep 0mm
\labelwidth 18mm
}%
\item[\rm (\theequation)\hfill]
#2
\end{list}%
}
\newcommand{\dyy}[2]{\dy{#1}{\raggedright$\dps#2$}}
\newcounter{stelling}
\newcommand{\thm}[2]{\refstepcounter{stelling}\vspace{4mm}\noindent{\bf Theorem \thestelling.}\label{#1}{\it #2}}
\newcounter{sectie}
\newcommand{\sect}[2]{\refstepcounter{sectie}
\section*{\boldmath \thesectie. #2}%
\label{#1}}
\newcommand{\sectz}[1]{\refstepcounter{sectie}
\section*{\boldmath \thesectie. #1}%
}
\newcommand{\pf}{\vspace{3mm}\noindent{\bf Proof.}\ }
\newcommand{\pfcl}{\vspace{3mm}\noindent{\em Proof.}\ }
\newcommand{\bx}{\hspace*{\fill} \hbox{\hskip 1pt \vrule width 4pt height 8pt depth 1.5pt \hskip 1pt}

\addvspace{4mm}}
\newcommand{\openbx}{\hspace*{\fill} $\Box$\\ \vspace{1mm}}
\newcommand{\rf}[1]{{\rm (\ref{#1})}}
\newcommand{\LABEL}[1]{\label{#1}}
\newcommand{\rank}{\text{\rm rank}}
\newcommand{\sgn}{\text{\rm sgn}}
\newcommand{\oF}{{\mathbb{F}}}

\newcommand{\oZ}{{\mathbb{Z}}}
\renewcommand{\loop}{\hspace*{-6pt}\raisebox{-.2\height}{\scalebox{0.06}{\includegraphics{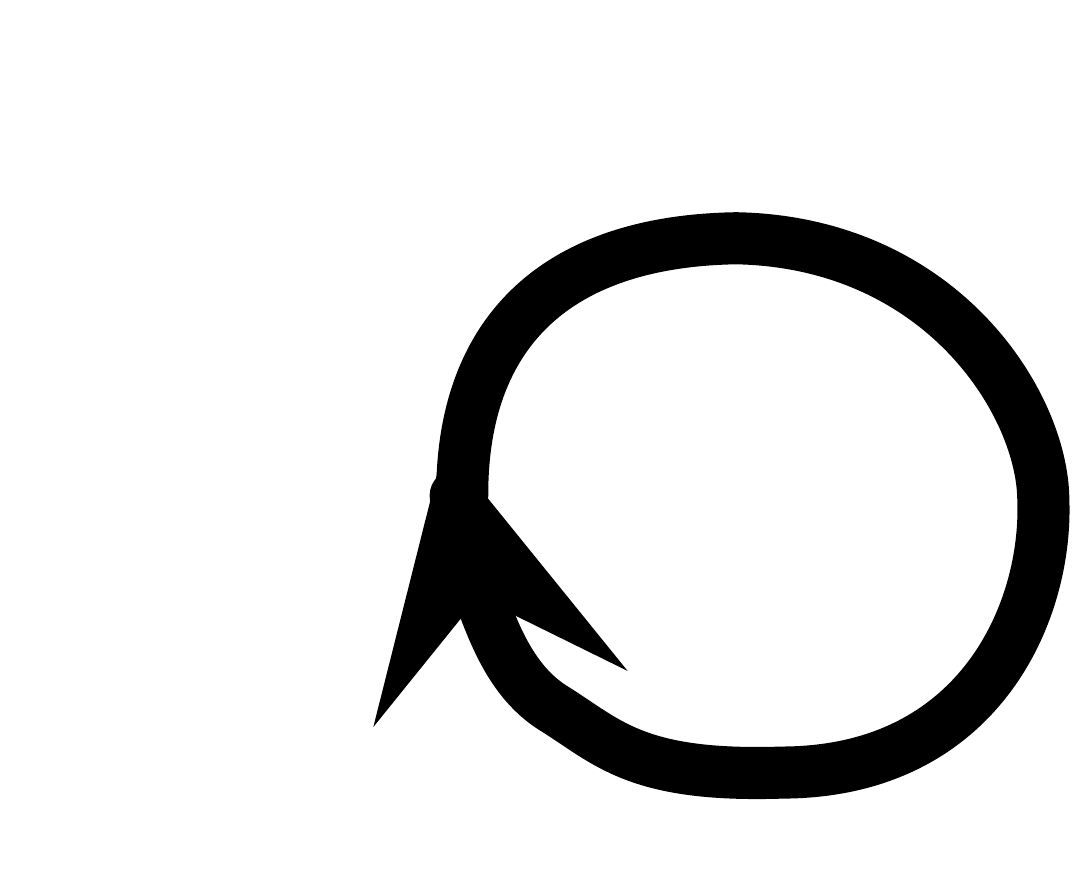}}}}
\newcommand{\Ker}{{\text{\rm Ker}}}
\renewcommand{\phi}{\varphi}
\begin{document}

\begin{center}
{\large\bf ON TRACES OF TENSOR REPRESENTATIONS OF DIAGRAMS

}
\vspace{3mm}
Alexander Schrijver \footnote{ University of Amsterdam and CWI, Amsterdam.
Mailing address:
Korteweg-de Vries Institute for Mathematics, University of Amsterdam,
P.O. Box 94248, 1090 GE Amsterdam, The Netherlands.
Email: lex@cwi.nl.
The research leading to these results has received funding from the European Research Council
under the European Union's Seventh Framework Programme (FP7/2007-2013) / ERC grant agreement
n$\mbox{}^{\circ}$ 339109.}

\end{center}

\medskip
\noindent
{\small
{\bf Abstract.}
Let $T$ be a set, of {\em types}, and let $\iota,o:T\to\oZ_+$.
A {\em $T$-diagram} is a locally ordered directed graph $G$
equipped with a function $\tau:V(G)\to T$ such that
each vertex $v$ of $G$ has indegree $\iota(\tau(v))$ and outdegree $o(\tau(v))$.
(A directed graph is {\em locally ordered} if
at each vertex $v$, linear orders of the edges entering $v$ and of the edges leaving $v$ are specified.)

Let $V$ be a finite-dimensional $\oF$-linear space, where $\oF$ is
an algebraically closed field of characteristic 0.
A function $R$ on $T$ assigning to each $t\in T$ a tensor
$R(t)\in V^{*\otimes \iota(t)}\otimes V^{\otimes o(t)}$ is called a {\em tensor representation} of $T$.
The {\em trace} (or {\em partition function}) of $R$
is the $\oF$-valued function $p_R$ on the collection of $T$-diagrams
obtained by `decorating' each vertex $v$ of a $T$-diagram $G$ with the tensor
$R(\tau(v))$, and contracting tensors along each edge of $G$, while respecting the order of the
edges entering $v$ and leaving $v$.
In this way we obtain a {\em tensor network}.

We characterize which functions on $T$-diagrams are traces, and show that each trace comes from a
unique `strongly nondegenerate' tensor representation.
The theorem applies to virtual knot diagrams, chord diagrams, and group representations.

}



\sectz{Introduction}

Our theorem characterizes traces of tensor networks, more precisely of tensor representations of
diagrams, which applies to knot diagrams, group representations, and algebras.
Tensor networks and their diagrammatical notation
root in work of Penrose [17], and were applied to
knot theory by Kauffman [9] and to Hopf algebra in `Kuperberg's notation' [12].
Other applications were found in areas like
quantum complexity
(cf.\ [1],
[8],
[15],
[19]),
statistical physics (cf.\ [7], [20]),
and
neural networks (cf.\ [16]).
(See Landsberg [13] for an in-depth survey of the geometry of tensors and its applications.)

\tussenkop{Types and $T$-diagrams}

Let $T$ be a (finite or infinite) set, of {\em types}, and let $\iota,o:T\to\oZ_+$ ($:=$ set
of nonnegative integers).
A {\em $T$-diagram} is a (finite) locally ordered directed graph $G$
equipped with a function $\tau:V(G)\to T$ such that
each vertex $v$ of $G$ has indegree $\iota(\tau(v))$ and outdegree $o(\tau(v))$.
Here a directed graph is {\em locally ordered} if
at each vertex $v$, a linear order of the edges entering $v$ and a linear order of the edges leaving $v$ are specified.
Loops and multiple edges are allowed.
Moreover, we allow the `vertexless directed loop' $\loop$ --- more precisely,
components of a $T$-diagram may be vertexless directed loops.

Let $\GG_T$ denote the collection of all $T$-diagrams.
If $T$ is clear from the context, we call a $T$-diagram just a {\em diagram}, and denote
$\GG_T$ by $\GG$.
The types can be visualized by small pictograms indicating the type of any vertex, as in the following examples.

\tussenkop{Examples}

{\em Virtual link diagrams.}
$T=\{\raisebox{-.2\height}{\scalebox{0.08}{\includegraphics{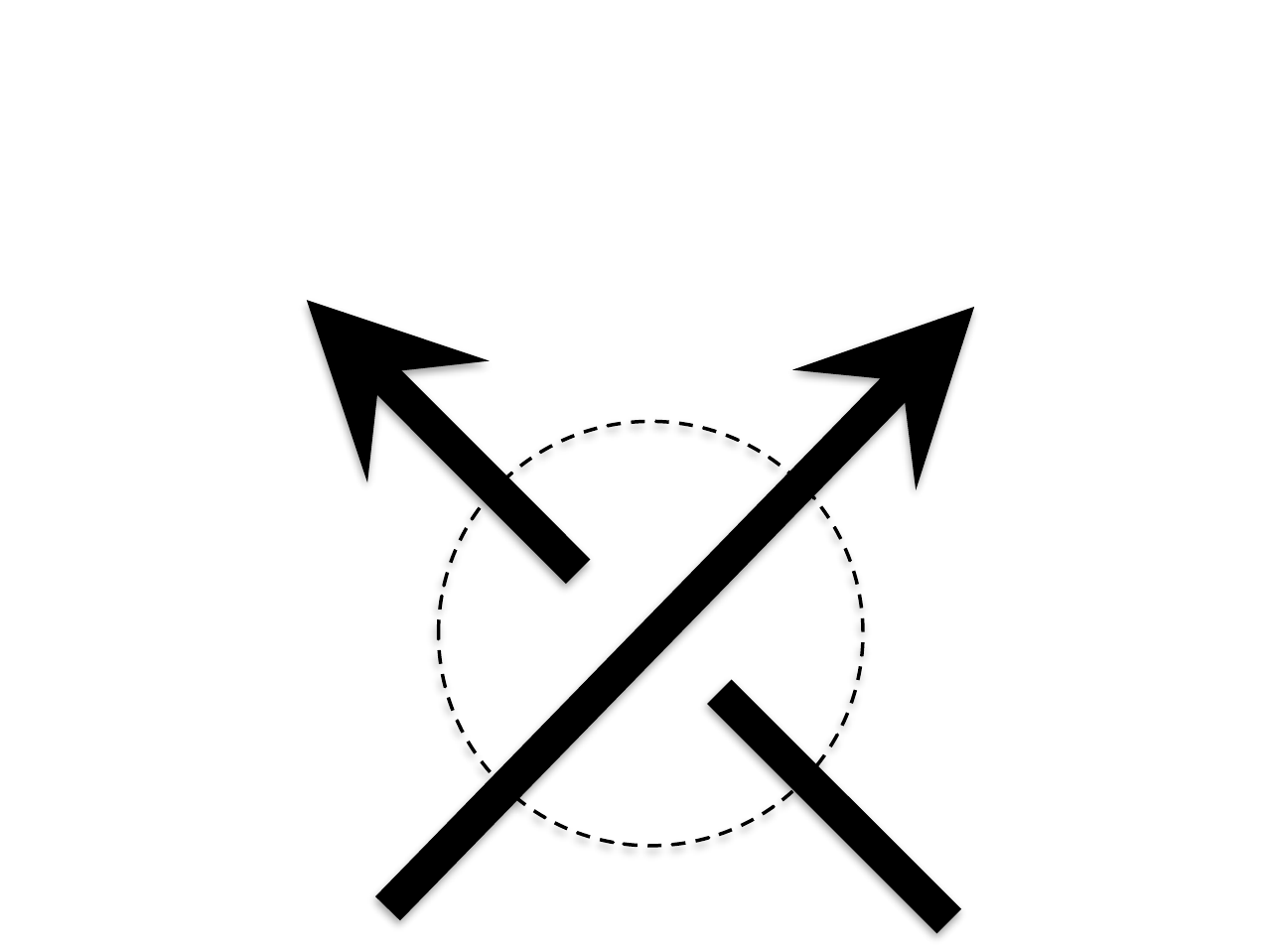}}},
\raisebox{-.2\height}{\scalebox{0.08}{\includegraphics{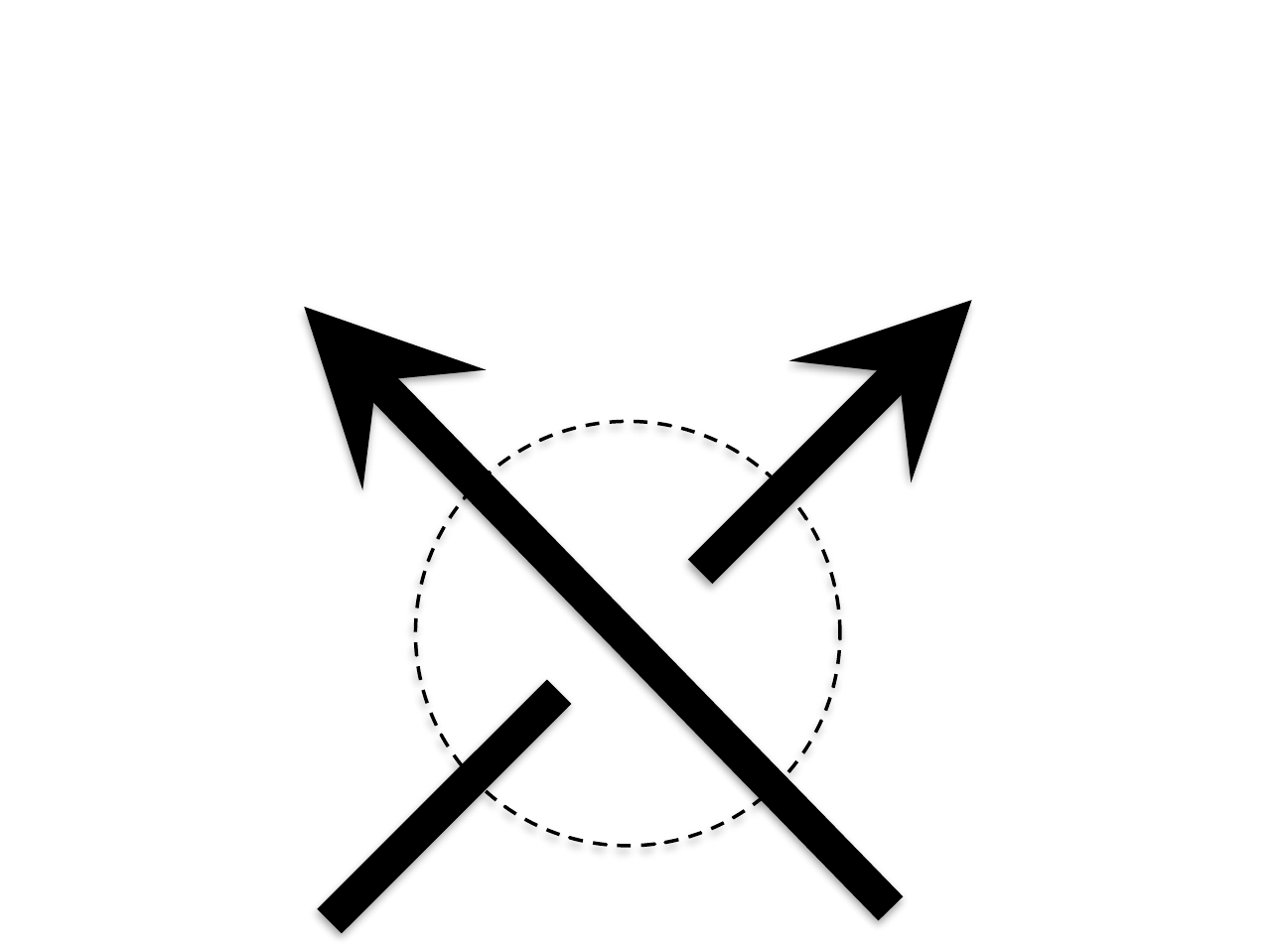}}}\}$.
So $|T|=2$ and $\iota(t)=o(t)=2$ for each $t\in T$.
(In pictures like this we assume the entering edges are ordered counter-clockwise and
the leaving edges are ordered clockwise.
We also will occasionally delete the grey circle indicating the vertex.)
Then the $T$-diagrams are the virtual link diagrams (cf.\ [10], [11], [14]).

\noindent
{\em Multiloop chord diagrams.}
$T=\{\raisebox{-.2\height}{\scalebox{0.1}{\includegraphics{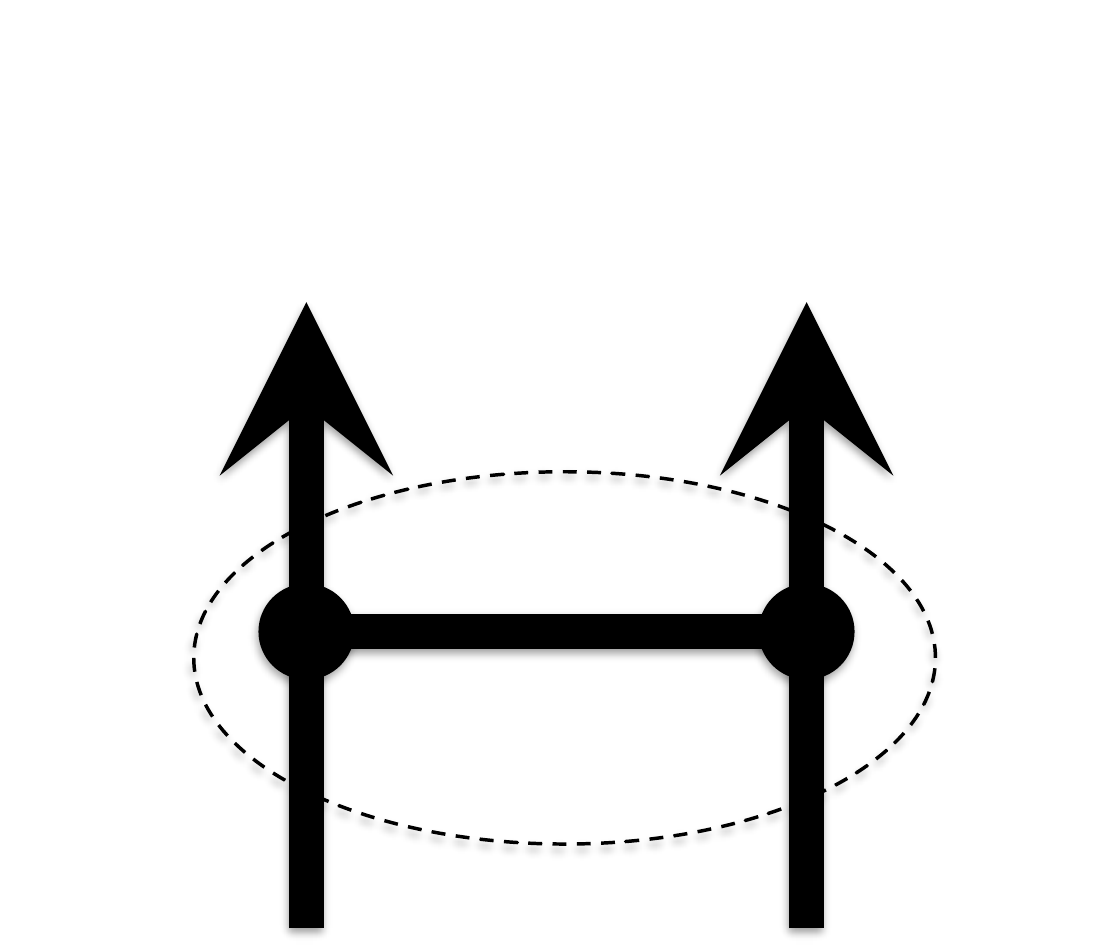}}}\}$, with
$\iota(\raisebox{-.2\height}{\scalebox{0.1}{\includegraphics{td_chord.pdf}}})
=
o(\raisebox{-.2\height}{\scalebox{0.1}{\includegraphics{td_chord.pdf}}})
=2$.
Then the $T$-diagrams are the multiloop chord diagrams, which play a key role in
the Vassiliev knot invariants (cf.\ [3]).
They can also be described as cubic graphs in which a set of disjoint oriented circuits (`Wilson loops')
covering all vertices is specified.
By contracting each Wilson loop to one point, the $T$-diagrams correspond to
graphs cellularly embedded on an oriented surface.

\smallskip
\noindent
{\em Groups.}
Let $\Gamma$ be a group, and let $T:=\Gamma$, with $\iota(t)=o(t)=1$ for each $t\in T$.
Then $T$-diagrams consist of disjoint directed cycles, with each vertex
typed by an element of $\Gamma$.

\smallskip
\noindent
{\em Algebra template.}
$T=\{
\hspace*{-8pt}\raisebox{-.3\height}{\scalebox{0.1}{\includegraphics{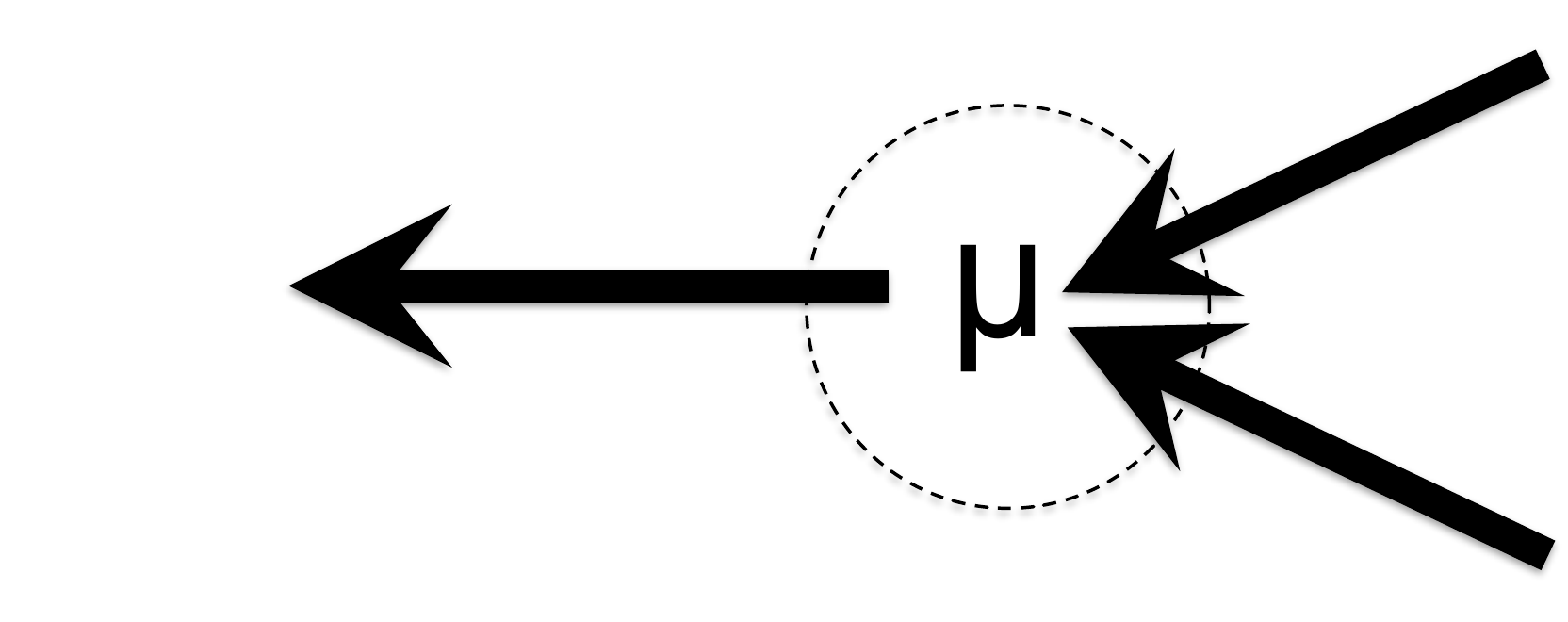}}}
,
\hspace*{-9pt}\raisebox{-.26\height}{\scalebox{0.1}{\includegraphics{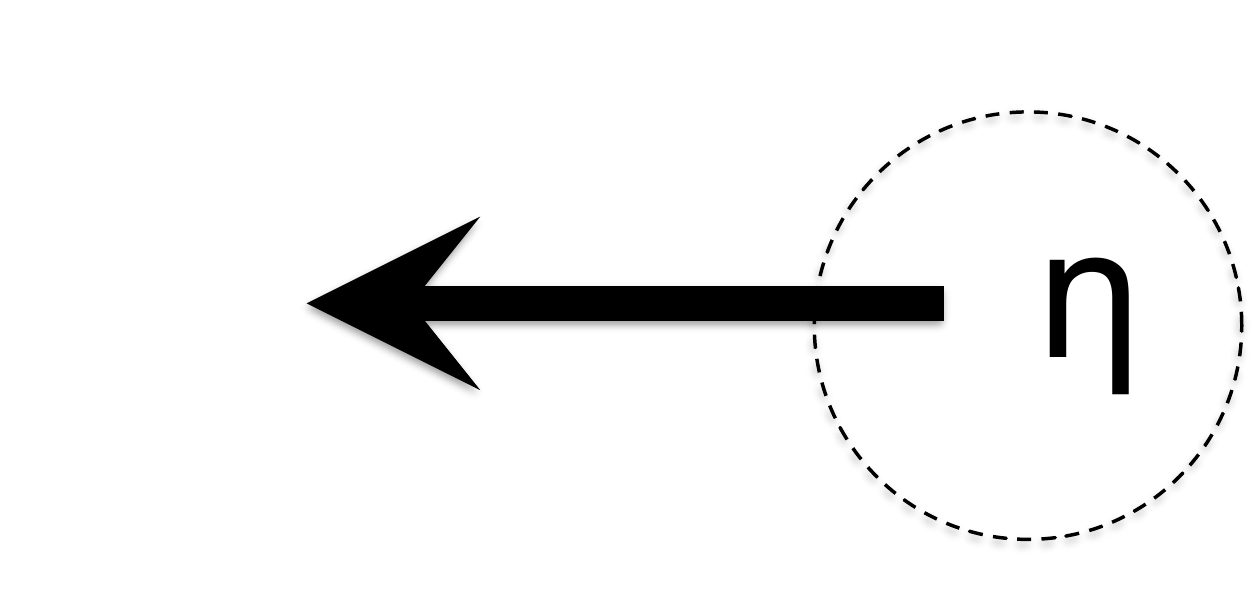}}}
\}$,
where the types represent the multiplication $\mu$ and the unit $\eta$, respectively.

\smallskip
\noindent
{\em Hopf algebra template.}
$T=\{
\hspace*{-8pt}\raisebox{-.3\height}{\scalebox{0.1}{\includegraphics{td_mu.pdf}}},
\hspace*{-9pt}\raisebox{-.26\height}{\scalebox{0.1}{\includegraphics{td_eta.pdf}}},
\hspace*{-7pt}\raisebox{-.36\height}{\scalebox{0.1}{\includegraphics{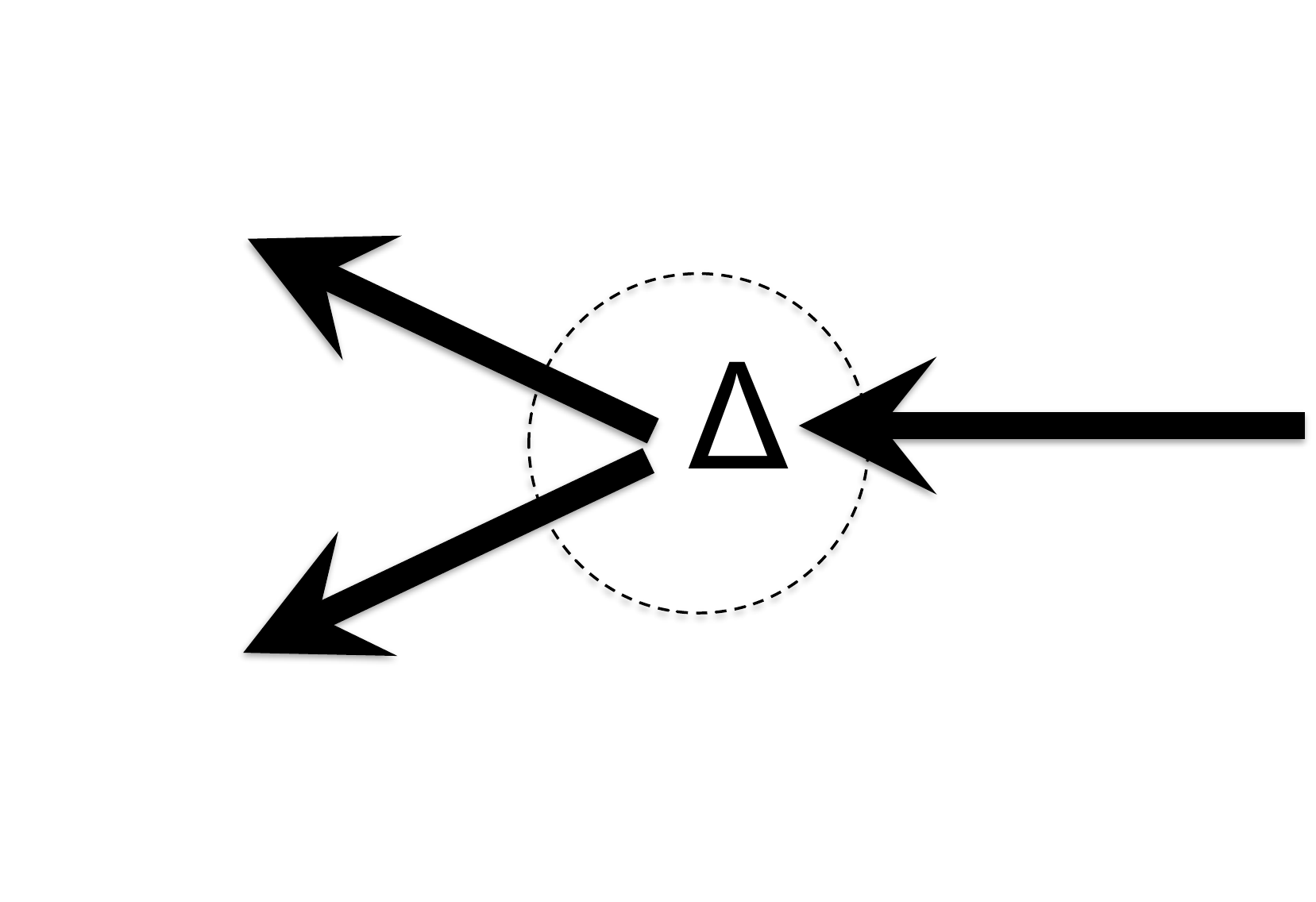}}},
\hspace*{-1pt}\raisebox{-.2\height}{\scalebox{0.1}{\includegraphics{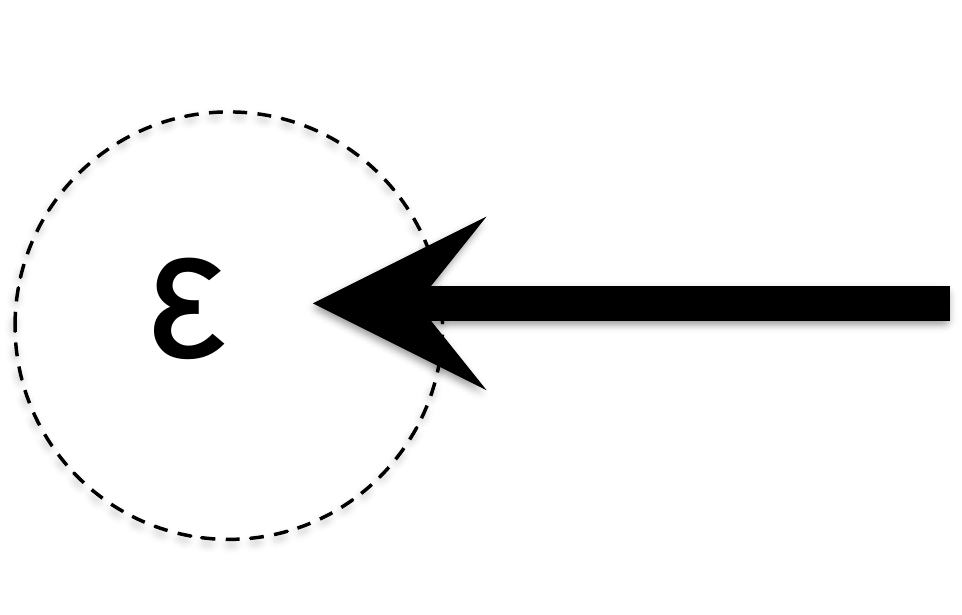}}},
\hspace*{-7pt}\raisebox{-.2\height}{\scalebox{0.1}{\includegraphics{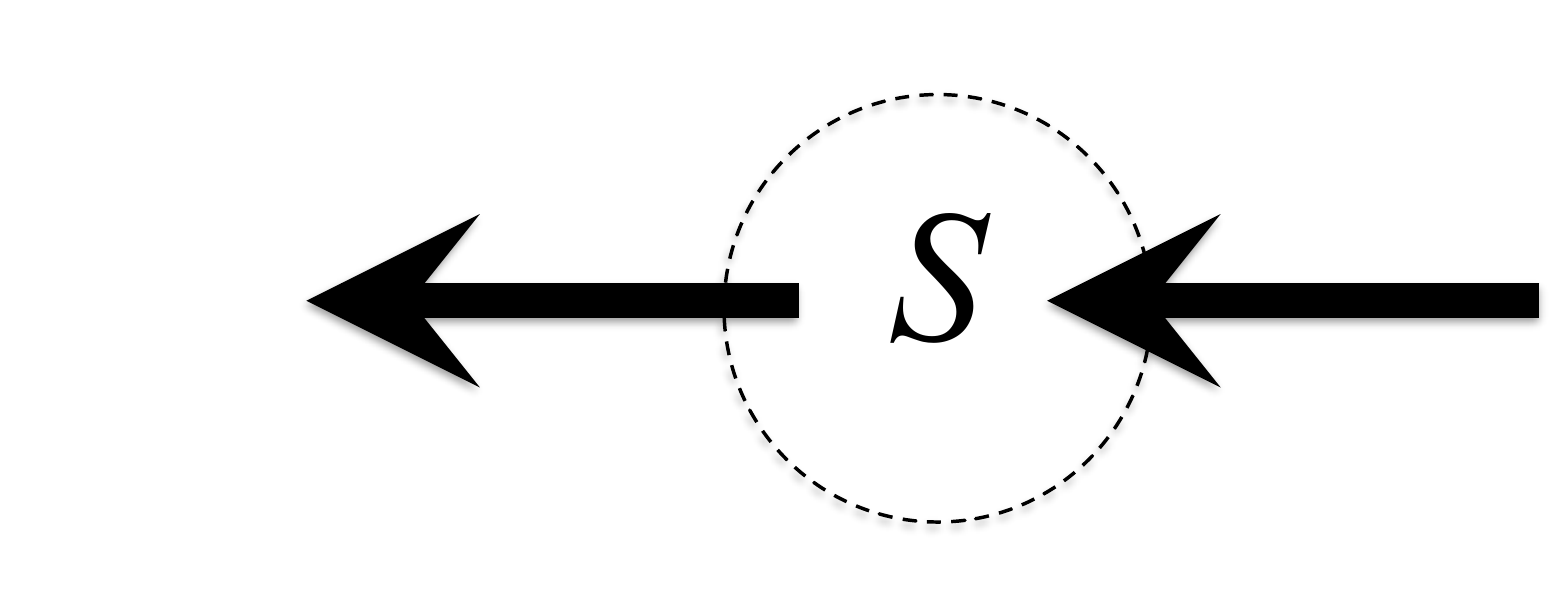}}}
\}$,
where
the types represent the multiplication $\mu$, the unit $\eta$, the comultiplication $\Delta$,
the counit $\varepsilon$, and the antipode $S$, respectively
(cf.\ Kuperberg [12]).

\smallskip
\noindent
{\em Directed graphs.}
$T:=\oZ_+^2$, with $\iota(k,l)=k$ and $o(k,l)=l$ for $(k,l)\in T$.

\tussenkop{Tensor representations and their traces}

Throughout this paper, fix an algebraically closed field $\oF$ of characteristic 0.
For any finite-dimensional $\oF$-linear space $V$, let as usual\footnote{We expect the two
different uses of $T$ as set of types and (non-italicized) in $\mbox{\rm T}(V)$ do not confuse.}
$$
\mbox{\rm T}(V):=\bigoplus_{k,l}(V^{*\otimes k}\otimes V^{\otimes l}).
$$
If $T$ is a set of types, call a function $R:T\to \mbox{\rm T}(V)$ a {\em tensor representation} of $T$
if $R(t)\in V^{*\otimes \iota(t)}\otimes V^{\otimes o(t)}$ for each $t\in T$.
We call $\dim(V)$ the {\em dimension} of $R$.
Let $\RR_T$ denote the collection of tensor representations $T\to \mbox{\rm T}(V)$.
($\RR_T$ depends on the linear space $V$, but we will use $\RR_T$ only when $V$ has been set.)

For a tensor representation $R:T\to \mbox{\rm T}(V)$, the {\em partition function}
or {\em trace} $p_R:\GG\to\oF$ of $R$ is defined as follows.
Roughly speaking, we `decorate' each vertex $v$ of a $T$-diagram $G$ with the tensor
$R(\tau(v))$, and contract tensors along each edge of $G$, consistent with the orders of the
edges entering $v$ and of those leaving $v$.
In this way we have a {\em tensor network}.

To give a more precise description of trace,
fix a basis $b_1,\ldots,b_n$ of $V$, with dual basis
$b^*_1,\ldots,b^*_n$.
Represent any element $x$ of
$V^{*\otimes k}\otimes V^{\otimes l}$ as a multi-dimensional array
$(x_{i_1,\ldots,i_k}^{j_1,\ldots,j_l})_{i_1,\ldots,i_k,j_1,\ldots,j_l=1}^n$, which are the coefficients
of $x$ when expressed in the basis $b^*_{i_1}\otimes\cdots\otimes b^*_{i_k}\otimes b_{j_1}\otimes\cdots\otimes b_{j_l}$
of $V^{*\otimes k}\otimes V^{\otimes l}$.
Set $[n]:=\{1,\ldots,n\}$.
Then
$$
p_R(G):=\sum_{\phi:E(G)\to [n]}\prod_{v\in V(G)}
R(\tau(v))_{\phi(\delta\inn(v))}^{\phi(\delta\uit(v))}.
$$
Here $\delta\inn(v)$ and $\delta\uit(v)$ are the ordered sets of edges entering $v$ and leaving $v$, respectively.
Moreover, for any ordered set $(e_1,\ldots,e_t)$ of edges, $\phi(e_1,\ldots,e_t):=(\phi(e_1),\ldots,\linebreak[1]\phi(e_t))$.

Note that $p_R(G)$ is independent of the chosen basis of $V$.
The function $\phi$ corresponds to a `state' or `edge coloring' of the `vertex model' $R$
of de la Harpe and Jones [7] (cf.\ [21]).

For $G\in\GG_T$, define $p(G):\RR_T\to\oF$ by $p(G)(R):=p_R(G)$.
Then $p(G)$ is $\GL(V)$-invariant, taking the natural action of $\GL(V)$ on $\RR_T$.
(We will use $p(G)$ only when $V$ has been set.)

\tussenkop{Webs = tangle diagrams}

To characterize which functions on the collection $\GG_T$ of $T$-diagrams are traces, we need the concept of tangle diagrams, also called webs.
When $T$ has been set,
for $k,l\in\oZ_+$, a {\em $k,l$-tangle diagram}, briefly a {\em $k,l$-web}, is a locally ordered directed graph $W$
equipped with injective functions $r:[k]\to V(W)$ and $s:[l]\to V(W)$ such that
$r(i)$ has outdegree 1 and
indegree 0 (for $i\in[k]$) and $s(j)$ has indegree 1 and outdegree 0 (for $j\in[l]$), and
equipped moreover with a function $\tau:V'(W)\to T$ such
that each vertex $v\in V'(W)$ has indegree $\iota(\tau(v))$ and outdegree $o(\tau(v))$,
where $V'(W):=V(W)\setminus(r([k])\cup s([l]))$.

The vertices in $r([k])$ are called the {\em roots} and the vertices in $s([l])$
are called the {\em sinks}.
For $i\in[k]$, $i$ is called the {\em label} of vertex $r(i)$, and for $j\in[l]$,
$j$ is called the {\em label} of vertex $s(j)$.
Again, loops and multiple edges are allowed, and components of $W$ may be the vertexless directe
loop \loop.
We call $W$ a {\em web} if it is a $k,l$-web for some $k,l$.
Let $\WW_{k,l}$ be the collection of all $k,l$-webs, and let $\WW$ be the collection
of all webs.
So $\WW_{0,0}=\GG$.
(We use this notation if $T$ has been set.)

By $\oF\GG$, $\oF\WW_{k,l}$, and $\oF\WW$ we denote the linear spaces of
formal $\oF$-linear combinations of elements of $\GG$, $\WW_{k,l}$, and $\WW$,
respectively.
Like in [5],
we call their elements {\em quantum diagrams}, {\em quantum $k,l$-webs}, and
{\em quantum webs}, respectively.
We extend any function on $\GG$, $\WW_{k,l}$, or $\WW$ to some linear space
linearly to a linear function on $\oF\GG$, $\oF\WW_{k,l}$, or $\oF\WW$.

For $G,H\in\GG$, let $G\cdot H$ be the disjoint union of $G$ and $H$.
More generally, if $W\in\WW_{k,l}$ and $X\in\WW_{l,k}$, let $W\cdot X$ be the diagram arising from
the disjoint union of $W$ and $X$ by,
for each $i\in[k]$, identifying the $i$-labeled root in $W$ with the $i$-labeled sink in $X$, and,
for each $j\in[l]$, identifying the $j$-labeled sink in $W$ with the $j$-labeled root in $X$.
After each identification, we ignore identified points
as vertex, joining the entering and leaving edge into one directed edge
(that is, \raisebox{-.28\height}{\scalebox{0.08}{\includegraphics{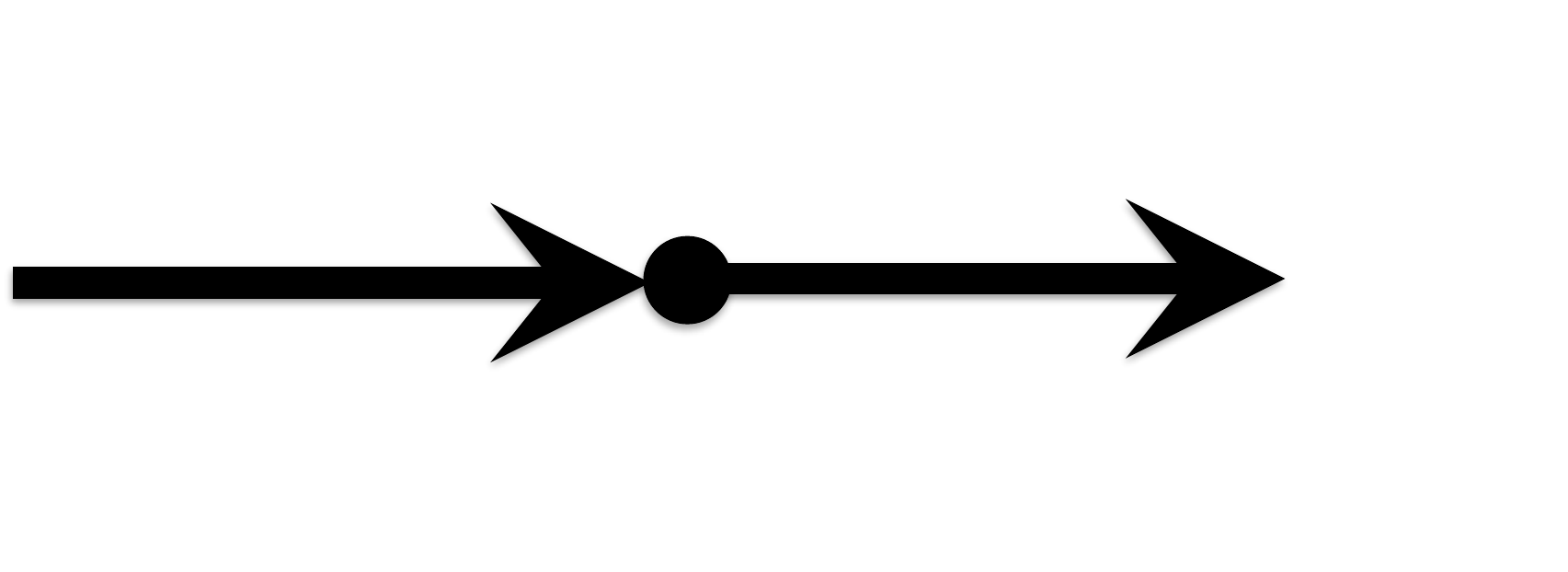}}}
becomes
\raisebox{-.28\height}{\scalebox{0.08}{\includegraphics{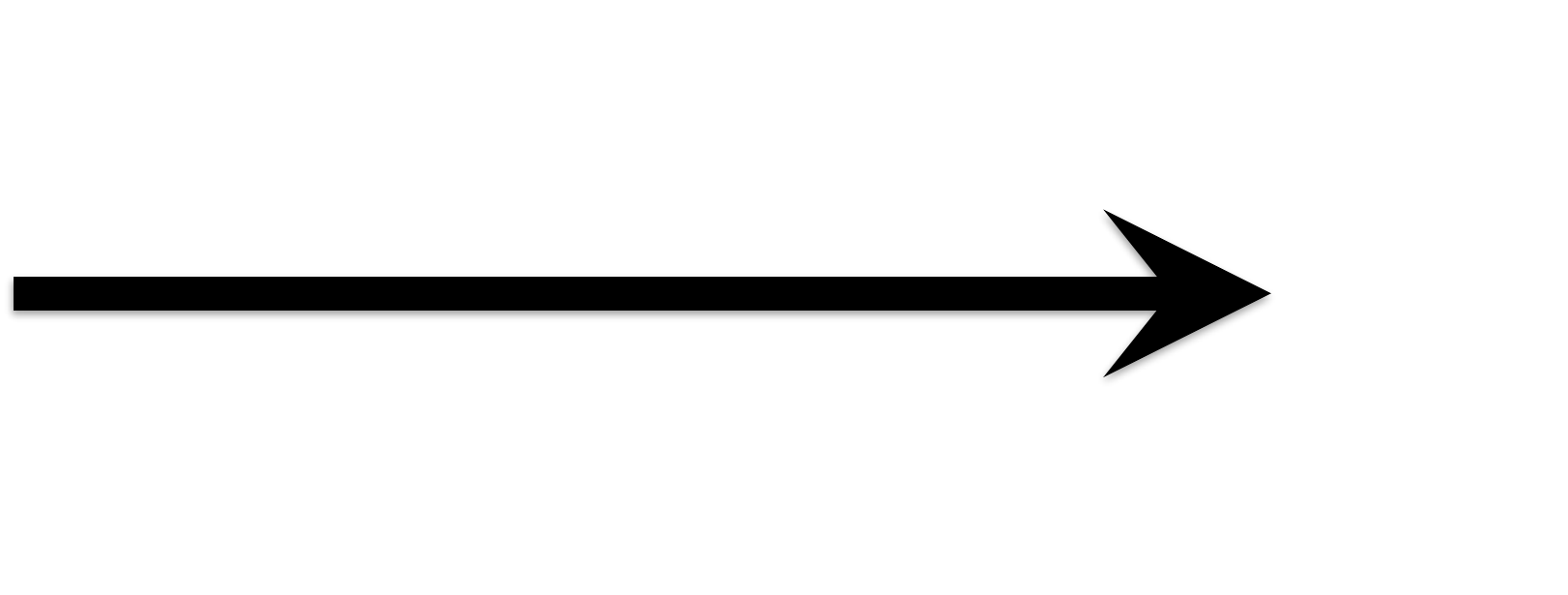}}}).
Note that this operation may introduce vertexless directed loops.
We extend this product $\cdot$ bilinearly to $\oF\WW\times\oF\WW\to\oF\WW$, setting
$W\cdot X:=0$ if $W\in\WW_{k,l}$ and $X\in\WW_{l',k'}$ with $(k,l)\neq (k',l')$.

Define, for each $k$, an element $\Delta_k\in\oF\WW_{k,k}$ as follows.
For $\pi\in S_k$ let $J_{\pi}$ be the $k,k$-web
consisting of $k$ disjoint directed edges $e_1,\ldots,e_k$, where
the tail of $e_i$ is labeled $i$ and its head is labeled $\pi(i)$,
for $i\in[k]$.
Then
$$
\Delta_k:=\sum_{\pi\in S_k}\sgn(\pi)J_{\pi}.
$$
Call $f:\GG\to\oF$ {\em multiplicative} if $f(\emptyset)=1$ and $f(G\cdot H)=f(G)f(H)$ for all $G,H\in\GG$.
Here $\emptyset$ is the diagram with no vertices and edges, and as before, $G\cdot H$ denotes the
disjoint union of $G$ and $H$.
We say that $f:\GG\to\oF$ {\em annihilates} a quantum web $\omega$ if
$f(\omega\cdot W)=0$ for each web $W$.

\thm{13me14a}{
Let $f:\GG_T\to\oF$.
Then there exists a tensor representation $R$ of $T$ of dimension $\leq n$ with $p_R=f$
if and only if $f$ is multiplicative and annihilates $\Delta_{n+1}$.
}

\bigskip
This theorem can be seen as a generalization of the following simple statement.
Let $\Gamma$ be any (finite or infinite) group.
Then a class function $\phi:\Gamma\to\oF$ is the character of some representation of
$\Gamma$ of dimension $\leq 2$ if and only if for all $a,b,c\in\Gamma$:%
$$
\phi(abc)+\phi(cba)+\phi(a)\phi(b)\phi(c)=
\phi(ab)\phi(c)+
\phi(ac)\phi(b)+
\phi(bc)\phi(a).
$$
(Similarly for higher dimensions.)

\tussenkop{The extended trace $\hat p_R$}

Generally, as the examples described above suggest,
we want to have a tensor representation that satisfies certain linear relations
between webs (for instance, `R-matrices' for the virtual link example).
Such relations can be described by a collection $Q$ of quantum webs.

Given a finite-dimensional $\oF$-linear space $V$ and a tensor representation $R:T\to \mbox{\rm T}(V)$,
we extend the trace function $p_R:\GG\to\oF$ to a function $\hat p_R:\WW\to \mbox{\rm T}(V)$ as follows.
Let $W$ be a $k,l$-web, with root function $r:[k]\to V(W)$ and sink function $s:[l]\to V(W)$.
Fix a basis $b_1,\ldots,b_n$ of $V$, with dual basis $b^*_1,\ldots,b^*_n$.
Then
\dyy{5ja15b}{
\hspace*{-2mm}
\hat p_R(W)
:=\sum_{\phi:E(W)\to[n]}\big(
\prod_{v\in V'(W)}
R(\tau(v))_{\phi(\delta\inn(v))}^{\phi(\delta\uit(v))}
\big)
\bigotimes_{i=1}^kb^*_{\phi(e'_i)}
\otimes
\bigotimes_{j=1}^lb_{\phi(e_j)},
}
where $V'(W):=V(W)\setminus(r([k])\cup s([l]))$,
and moreover, for $i\in[k]$, $e'_i$ is the edge leaving $r(i)$,
and, for $j\in[l]$, $e_j$ is the edge entering $s(j)$.

Again, $\hat p_R(W)$ is independent of the chosen basis of $V$.
Also, $p_R:=\hat p_R|\GG$.
We set $\hat p(W)(R):=\hat p_R(W)$ for $R\in\RR_T$ and web $W$.
Then for each $W$, $\hat p(W)$ is $\GL(V)$-invariant.
Finally, by letting $\cdot$ (also) to be the standard bilinear form on $\mbox{\rm T}(V)$, we have for webs $W$ and $X$:
\dyy{5ja15a}{
\hat p_R(W)\cdot\hat p_R(X)=p_R(W\cdot X).
}
Hence for any set $Q$ of quantum webs, $\hat p_R(Q)=0$ implies that $p_R$ annihilates $Q$
(meaning that it annihilates each $\omega\in Q$) ---
but not conversely.
(An easy example is $T:=\{a\}$ with $\iota(a)=o(a)=1$, $V=\oF^2$, and
$R:=
{\tiny
\left(
\begin{array}{cc}
1&1\\
0&1
\end{array}
\right)}$,
$Q=\{
\hspace*{-5pt}\raisebox{-.38\height}{\scalebox{0.12}{\includegraphics{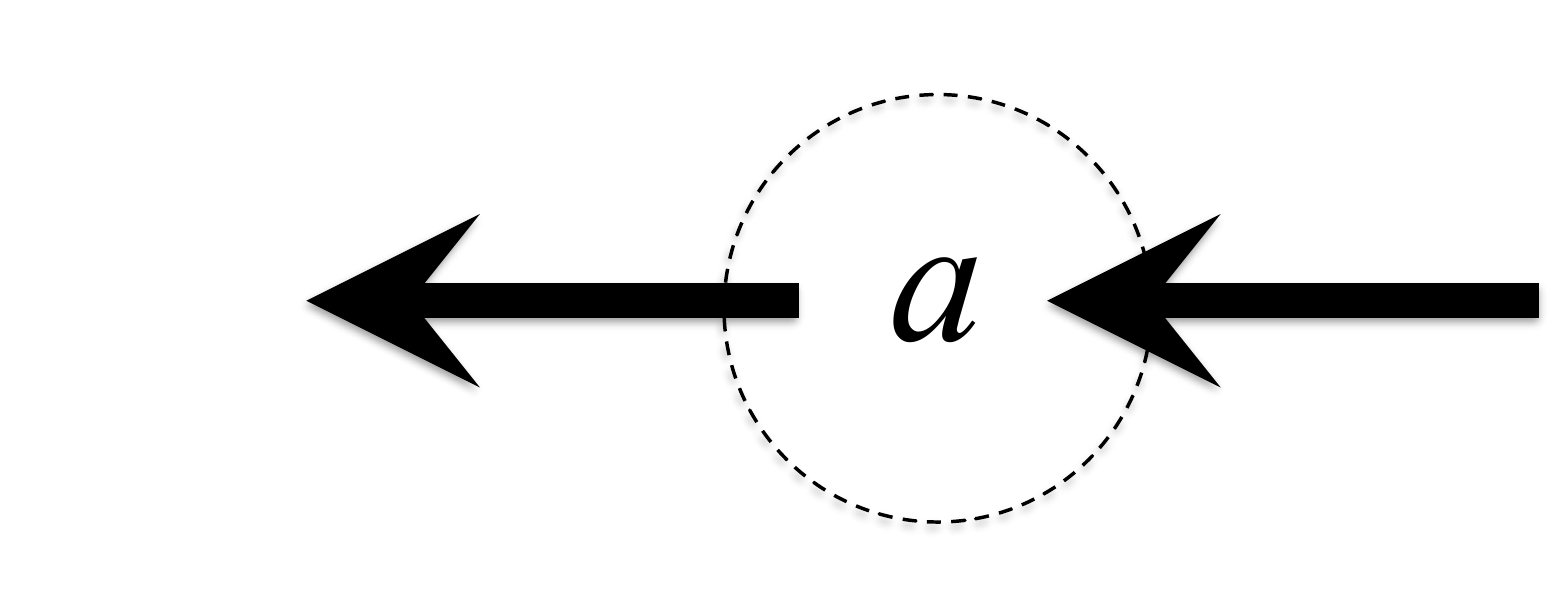}}}
-
\hspace*{-9pt}\raisebox{-.357\height}{\scalebox{0.12}{\includegraphics{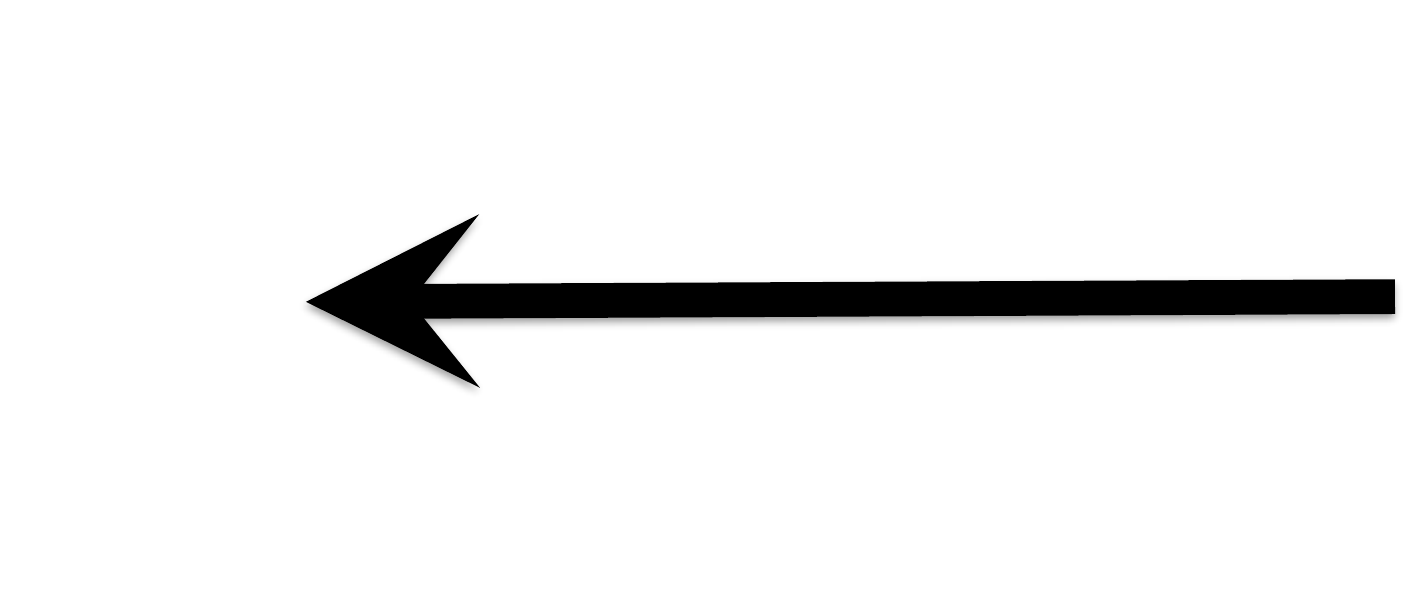}}}
\}$.)

However, as we will see, if $p_R$ annihilates $Q$,
then there exists $R'$ such that $p_{R'}=p_R$ and $\hat p_{R'}(\omega)=0$ for all $\omega\in Q$.
So we could take for $Q$ the collection of {\em all} quantum webs annihilated by $p_R$.

\tussenkop{Examples (continued)}

{\em Virtual link diagrams.}
The following set of quantum webs correspond to the Reidemeister moves:
$$
Q:=
\{
\hspace*{-7pt}\raisebox{-.3\height}{\scalebox{0.1}{\includegraphics{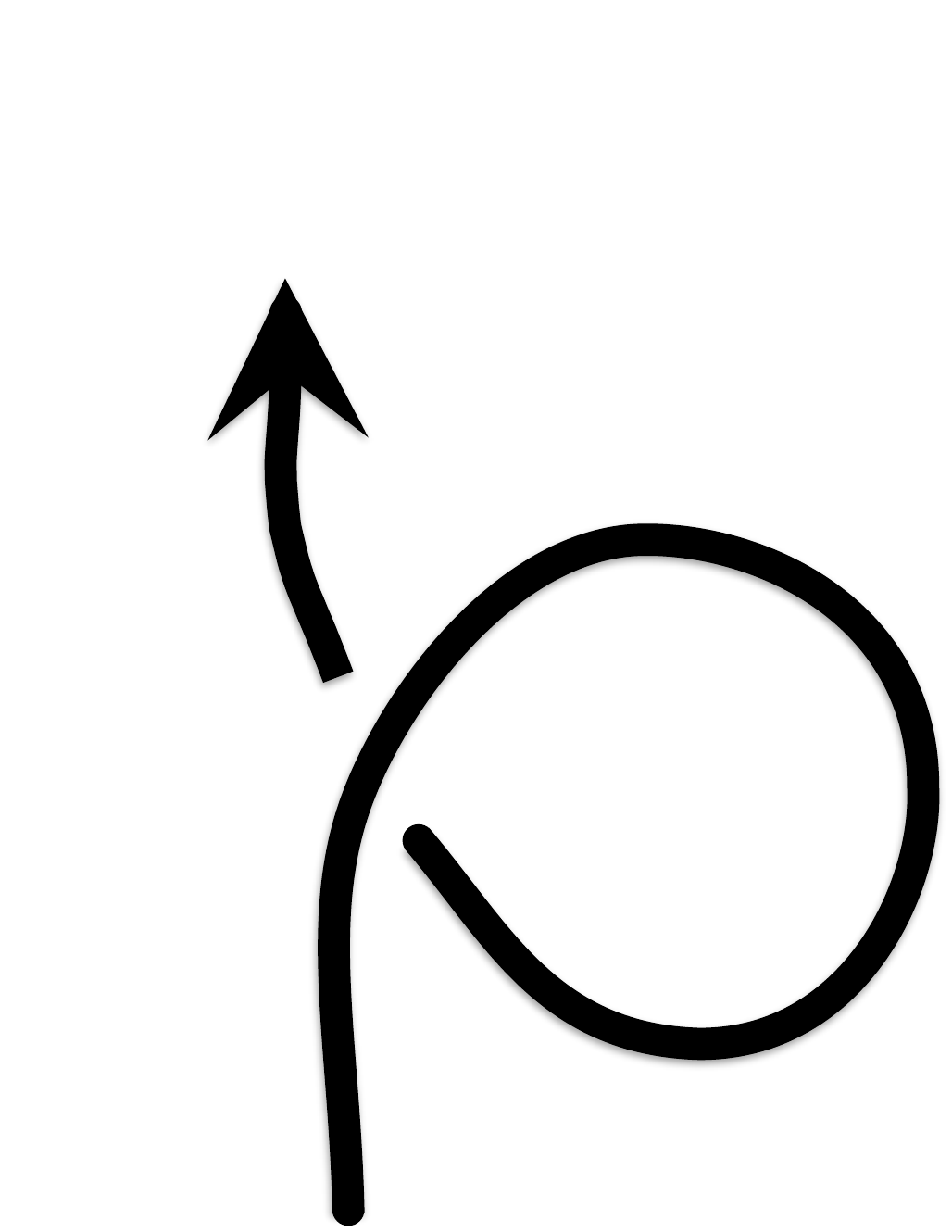}}}
-
\hspace*{-9pt}\raisebox{-.3\height}{\scalebox{0.1}{\includegraphics{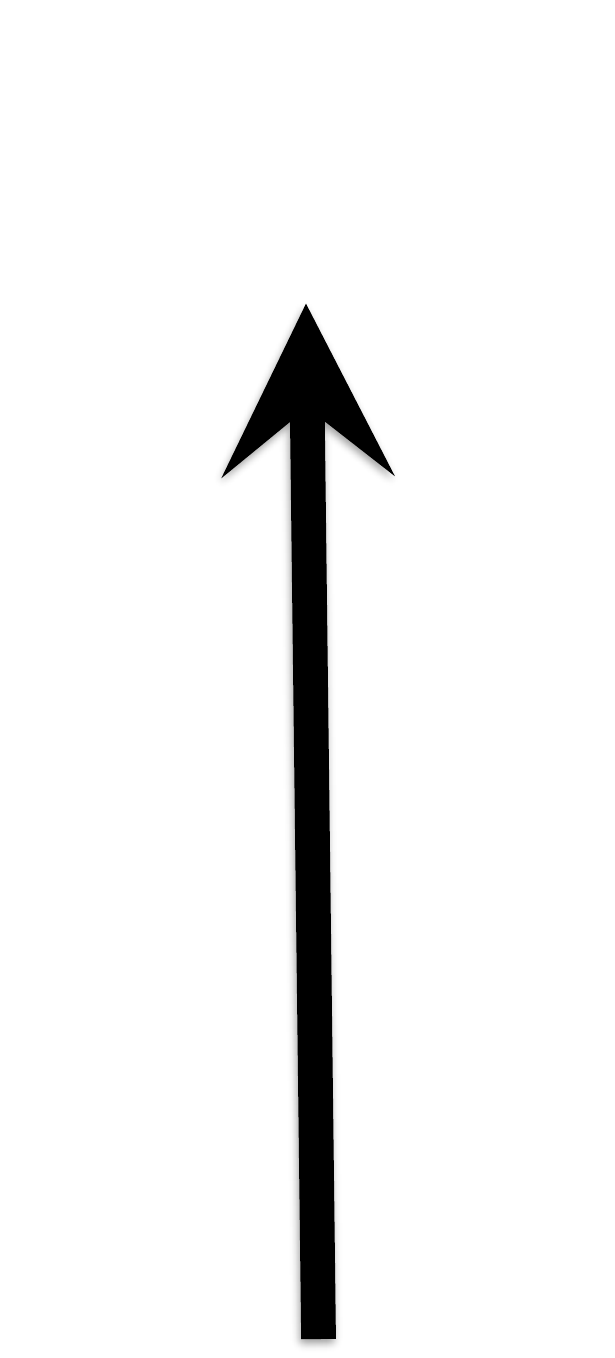}}},
\hspace*{-9pt}\raisebox{-.3\height}{\scalebox{0.1}{\includegraphics{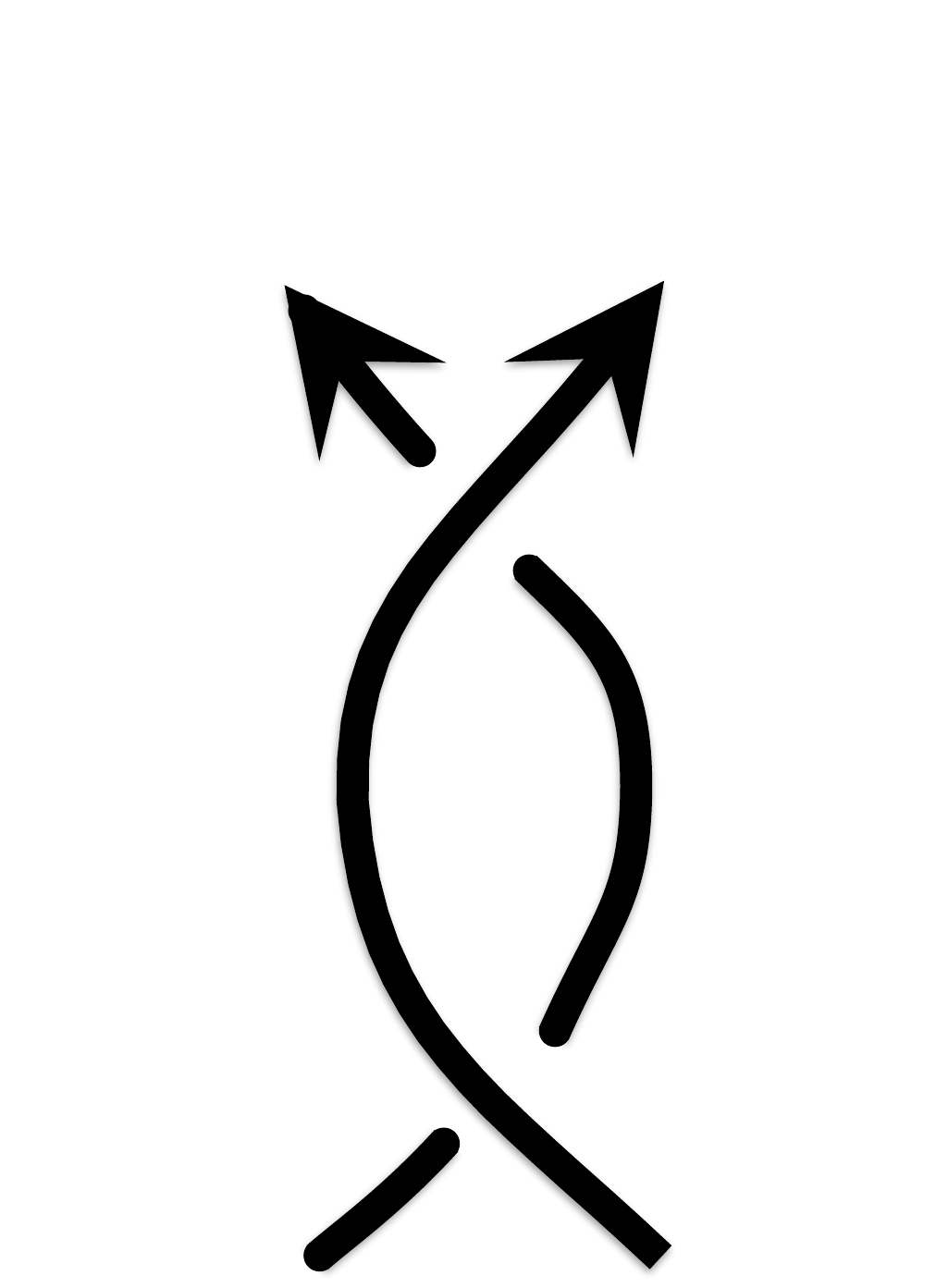}}}
-
\hspace*{-9pt}\raisebox{-.3\height}{\scalebox{0.1}{\includegraphics{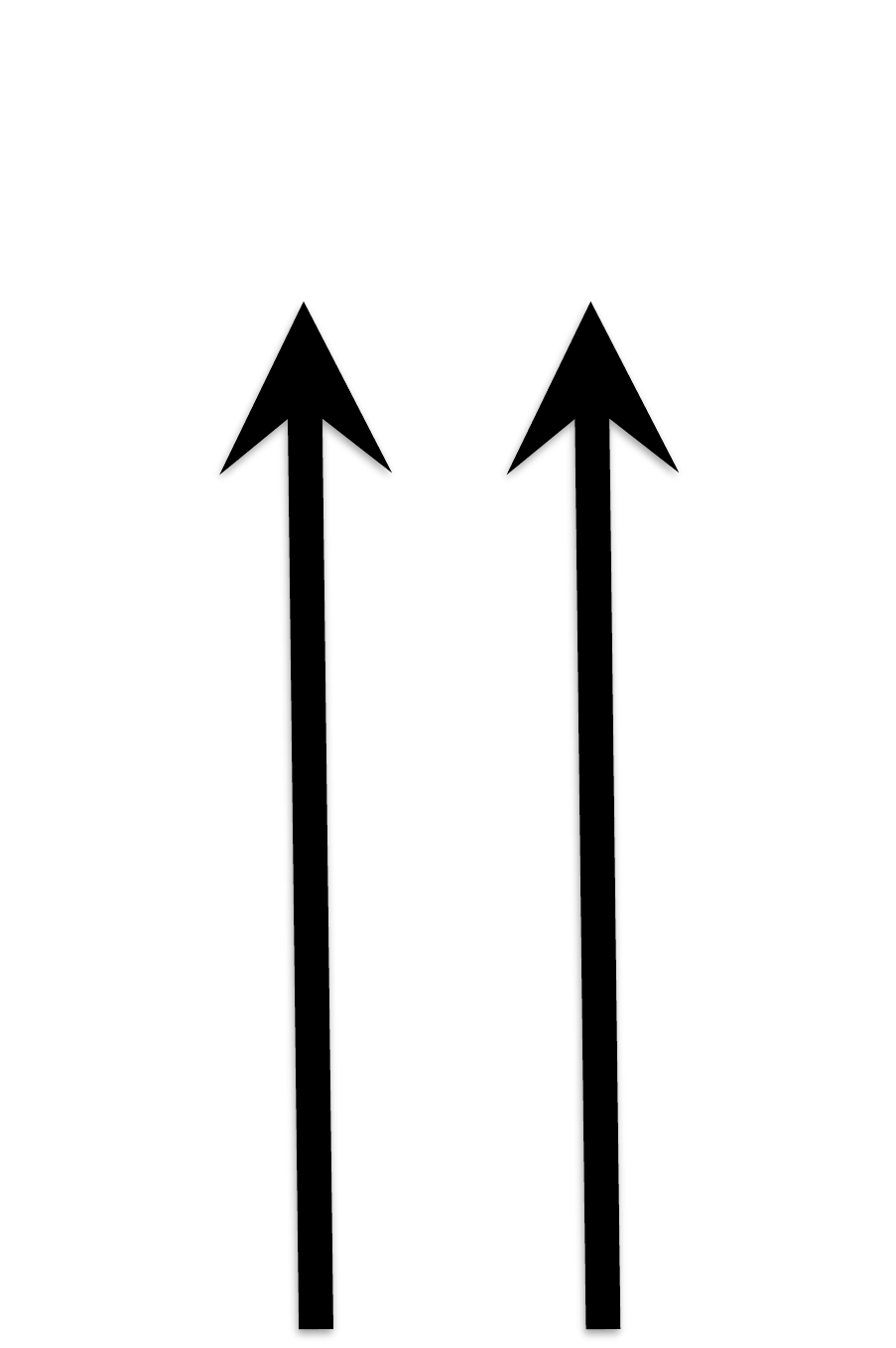}}},
\hspace*{-9pt}\raisebox{-.4\height}{\scalebox{0.1}{\includegraphics{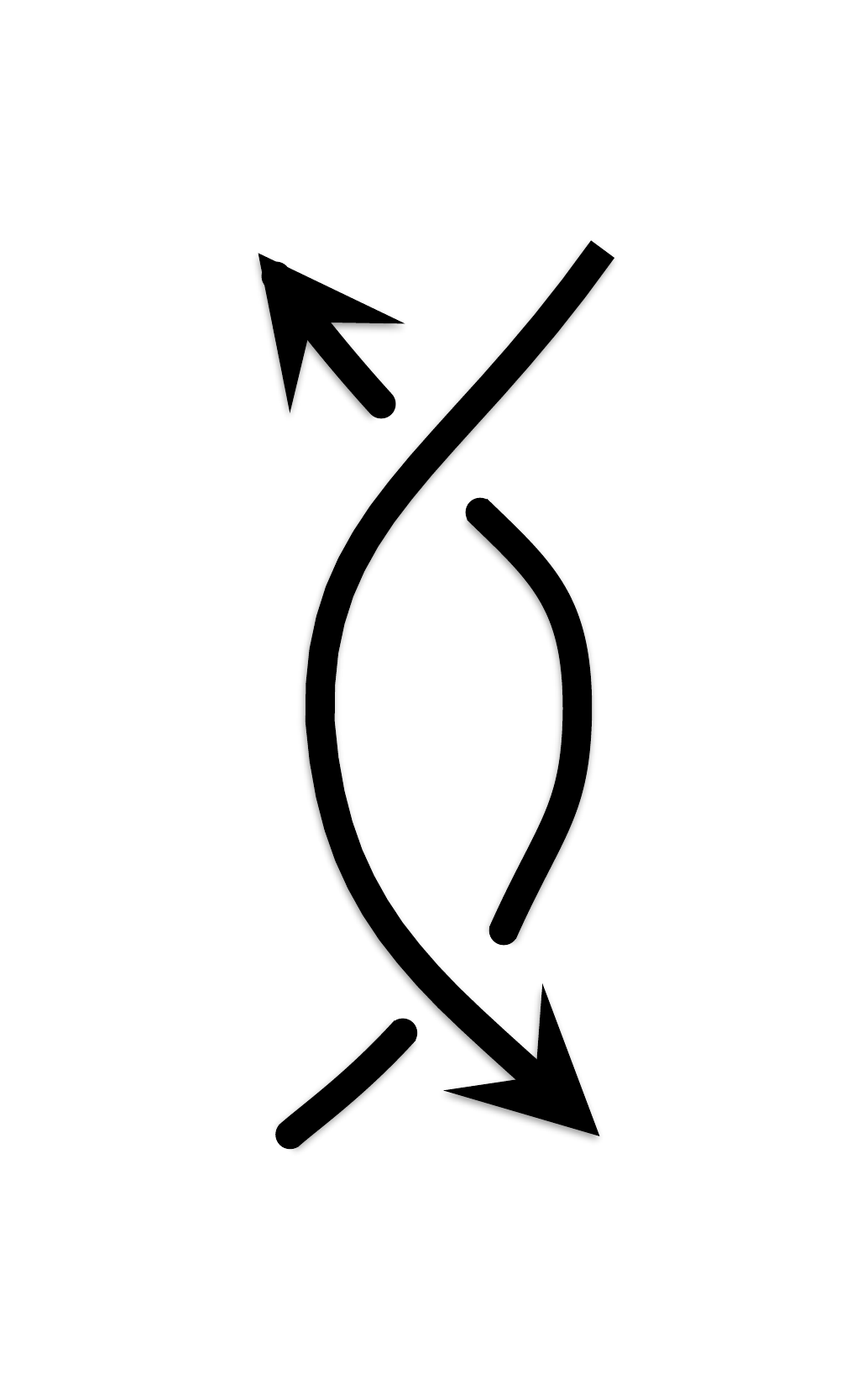}}}
-
\hspace*{-9pt}\raisebox{-.4\height}{\scalebox{0.1}{\includegraphics{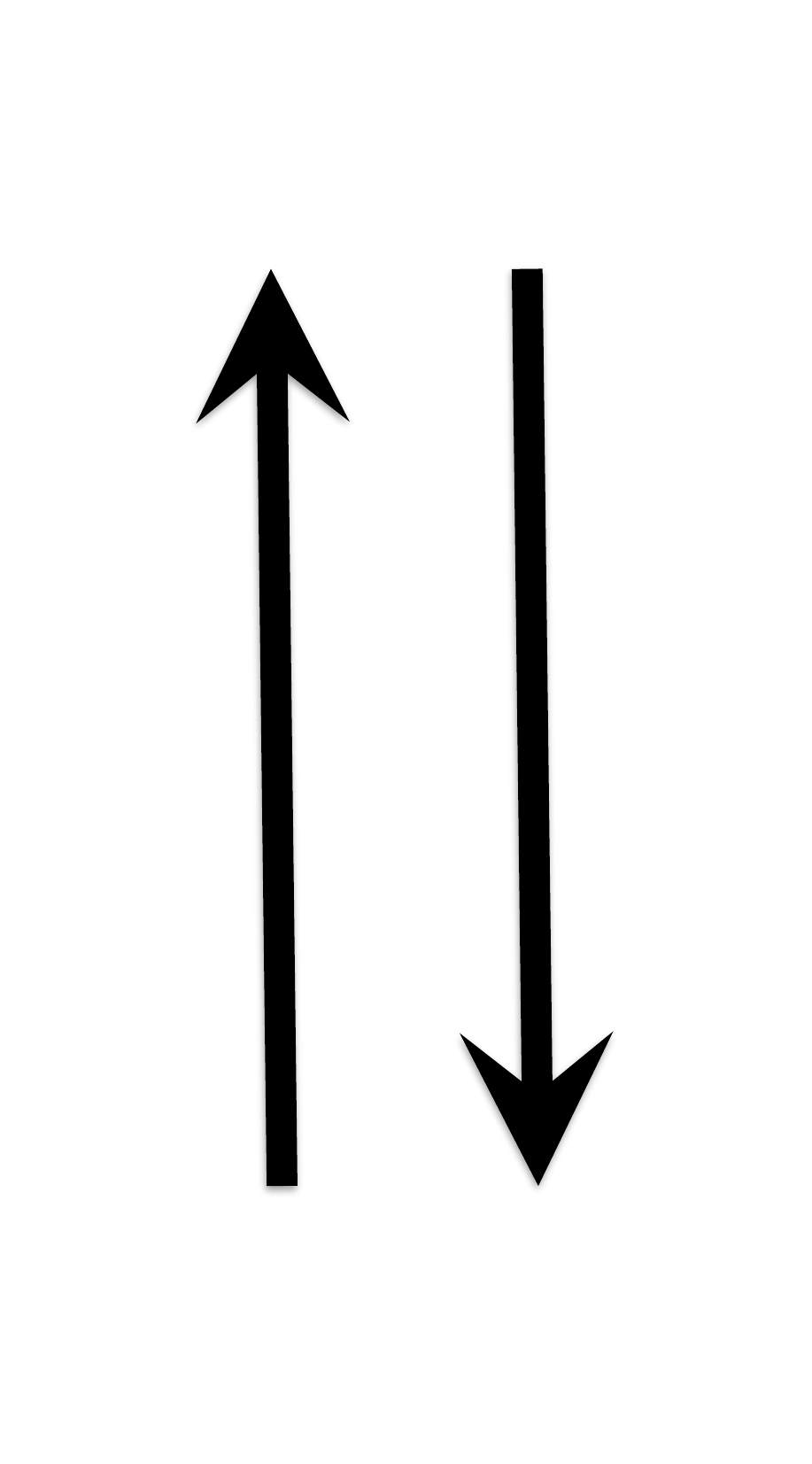}}},
\hspace*{-9pt}\raisebox{-.3\height}{\scalebox{0.1}{\includegraphics{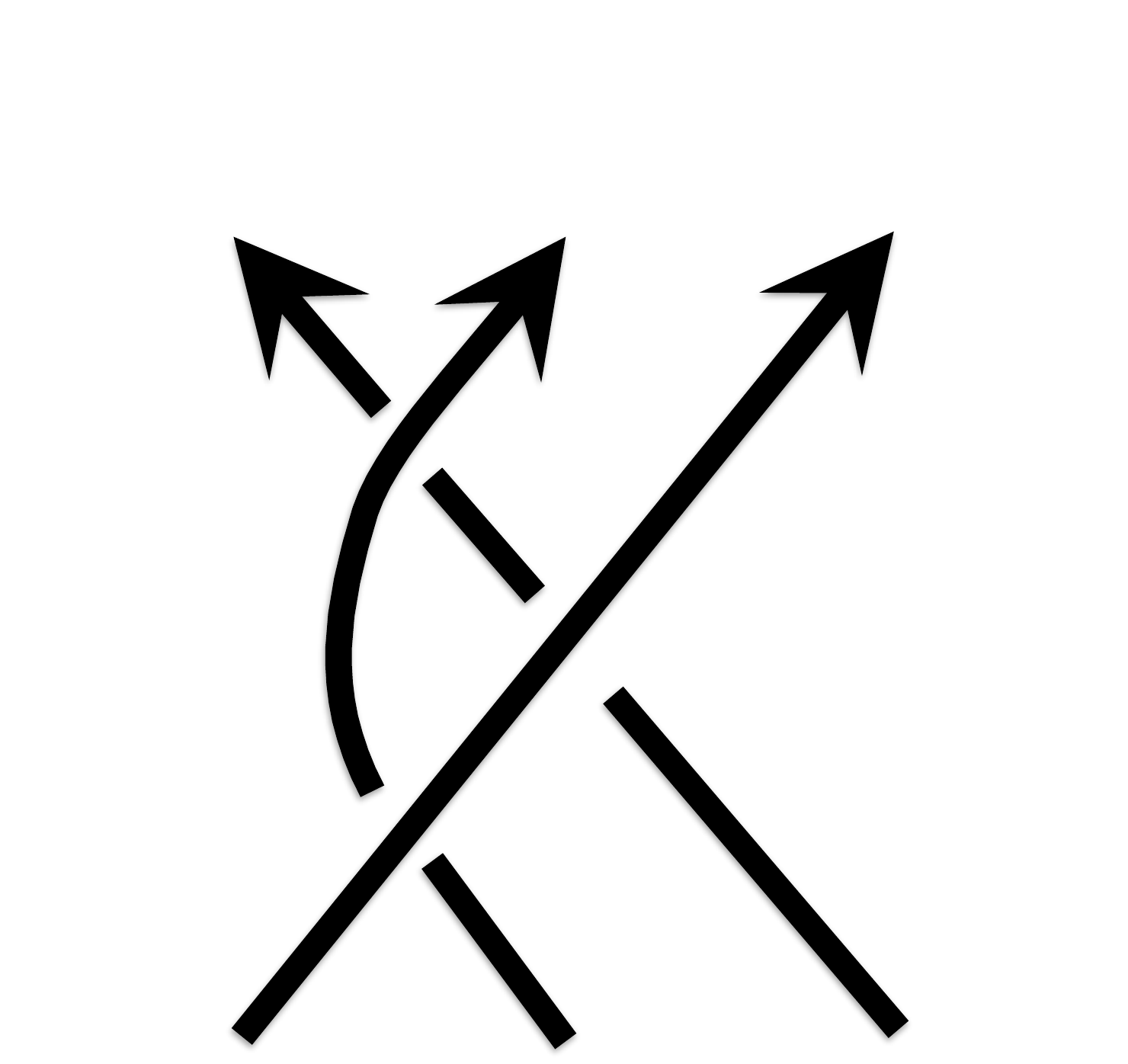}}}
-
\hspace*{-9pt}\raisebox{-.3\height}{\scalebox{0.1}{\includegraphics{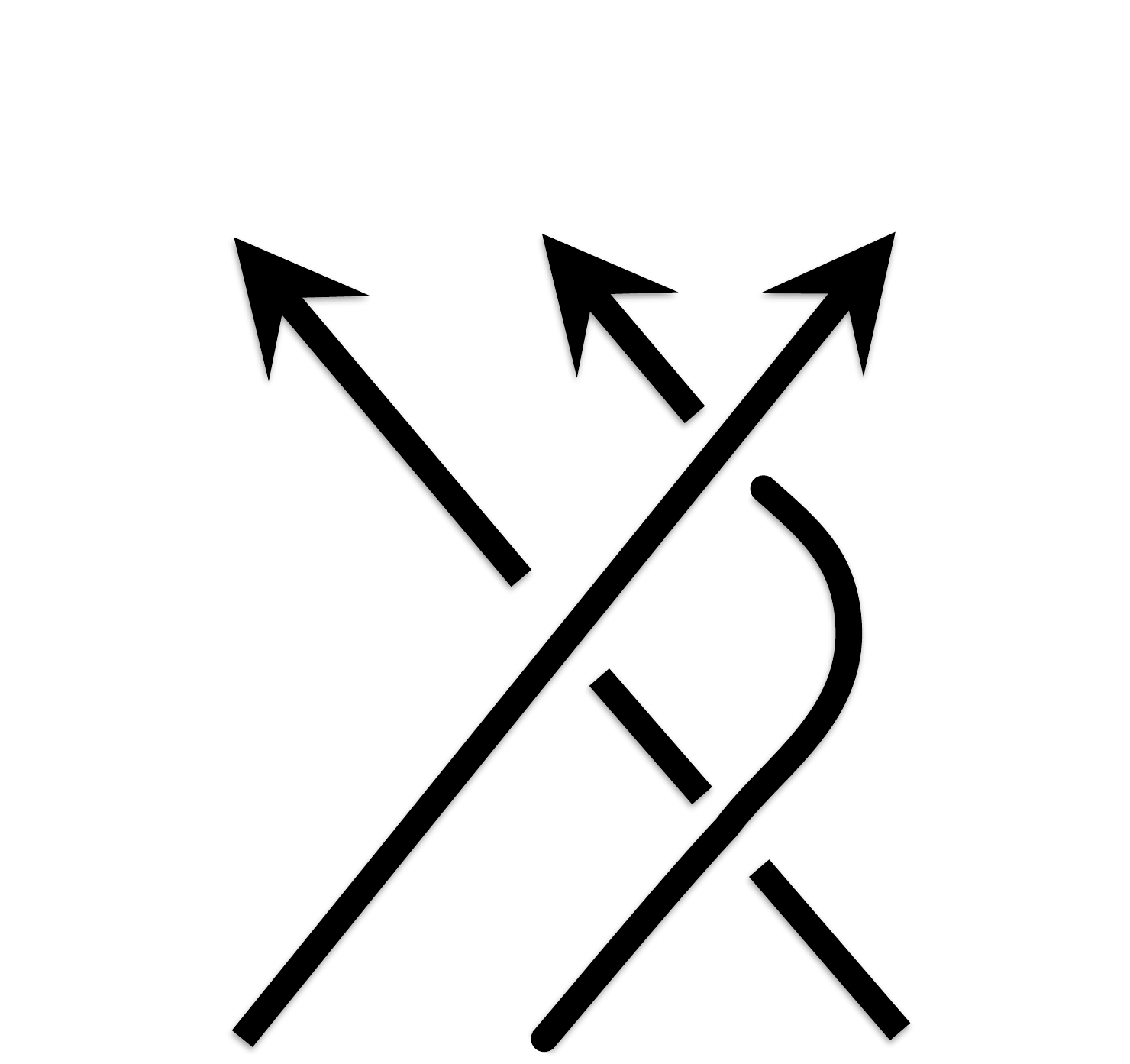}}}
\}.
$$
(In pictures representing quantum webs like this we assume that the roots and sinks in the
different webs occurring in the quantum web are labeled consistently suggested by their
position in the pictures.
The precise numbering of the roots and sinks is irrelevant, as long it is consistent over all webs
occurring in the quantum web.)
Then the functions $f:\GG\to\oF$ annihilating $Q$ are the virtual link invariants
(that is, invariant under Reidemeister moves).
Moreover, $\hat p_R(Q)=0$ if and only if $R$ is an `R-matrix'.

\smallskip
\noindent
{\em Multiloop chord diagrams.}
To describe the `undirectedness' of the chords and the `4T-relations', set
$$
Q:=
\{
\hspace*{-3pt}\raisebox{-.3\height}{\scalebox{0.1}{\includegraphics{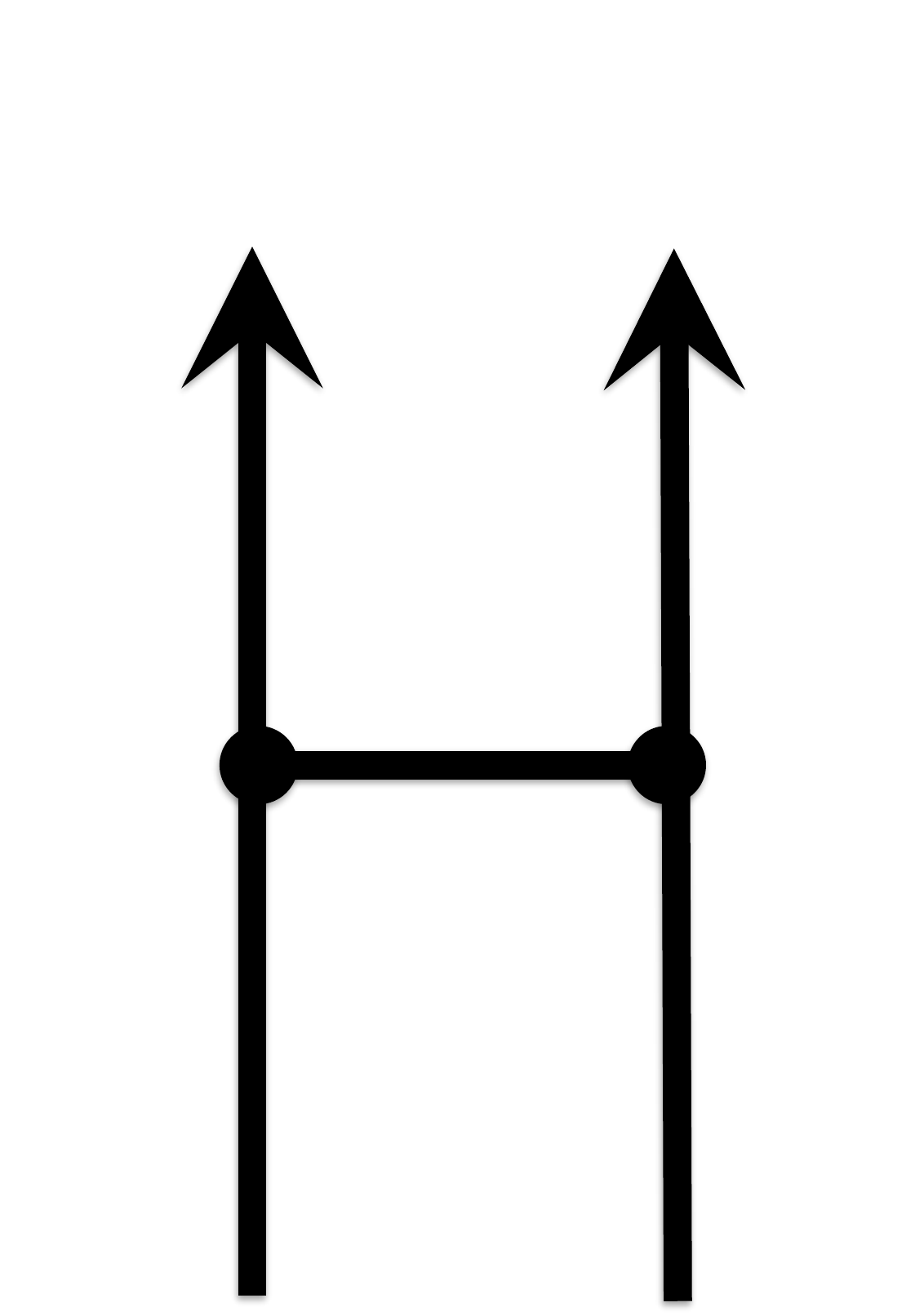}}}
-
\hspace*{-3pt}\raisebox{-.3\height}{\scalebox{0.1}{\includegraphics{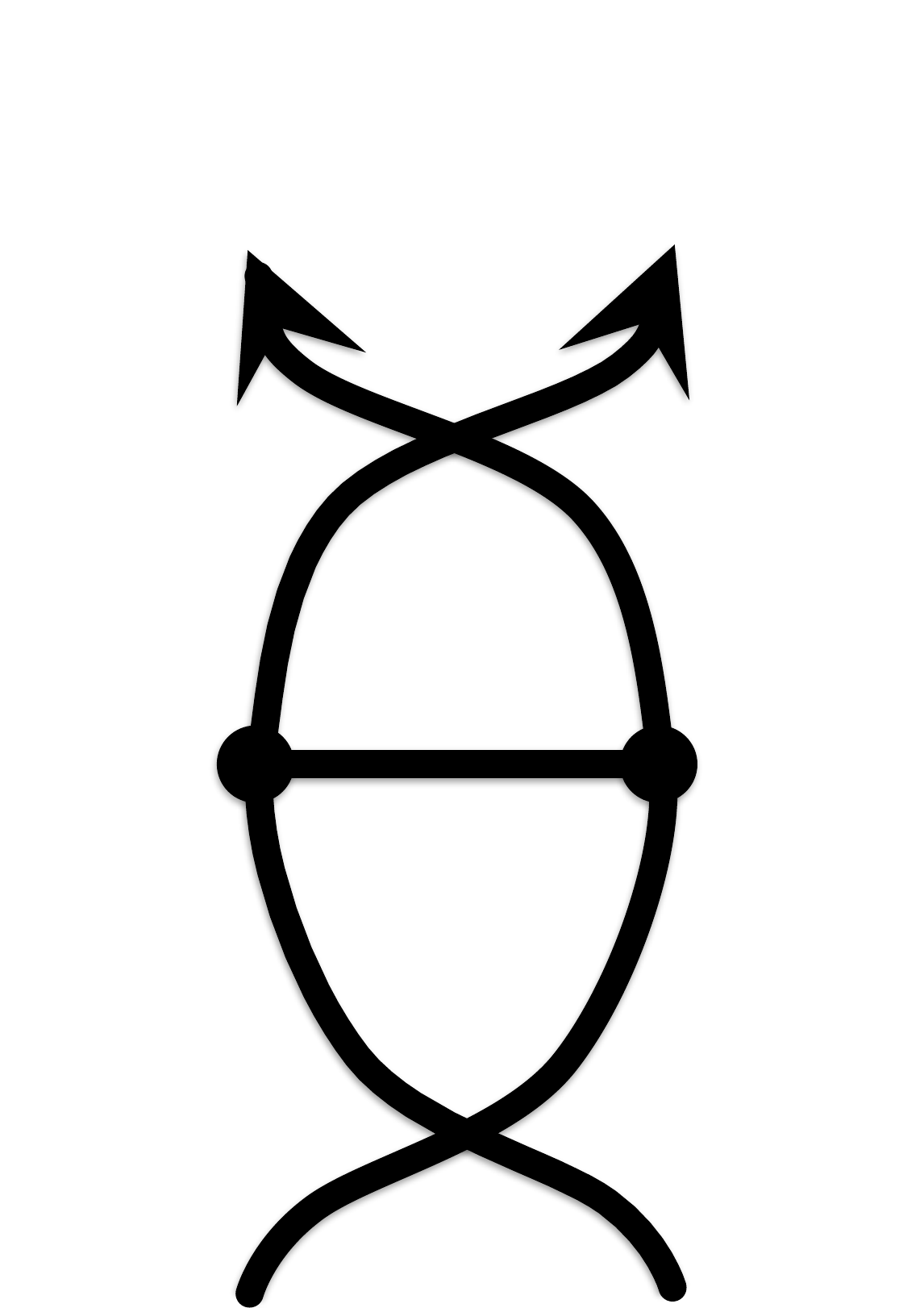}}}
,
\hspace*{-3pt}\raisebox{-.3\height}{\scalebox{0.1}{\includegraphics{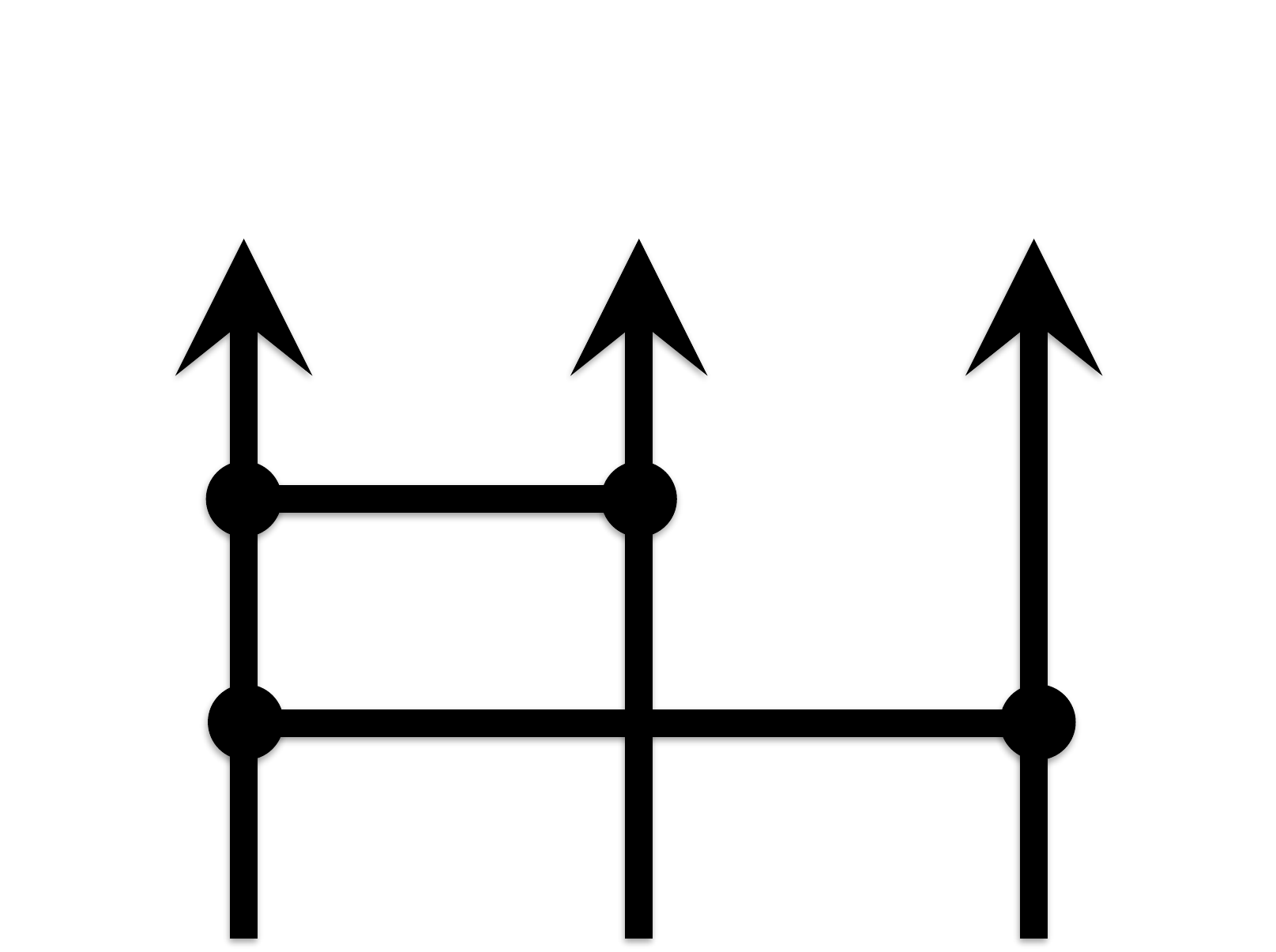}}}
+
\hspace*{-3pt}\raisebox{-.3\height}{\scalebox{0.1}{\includegraphics{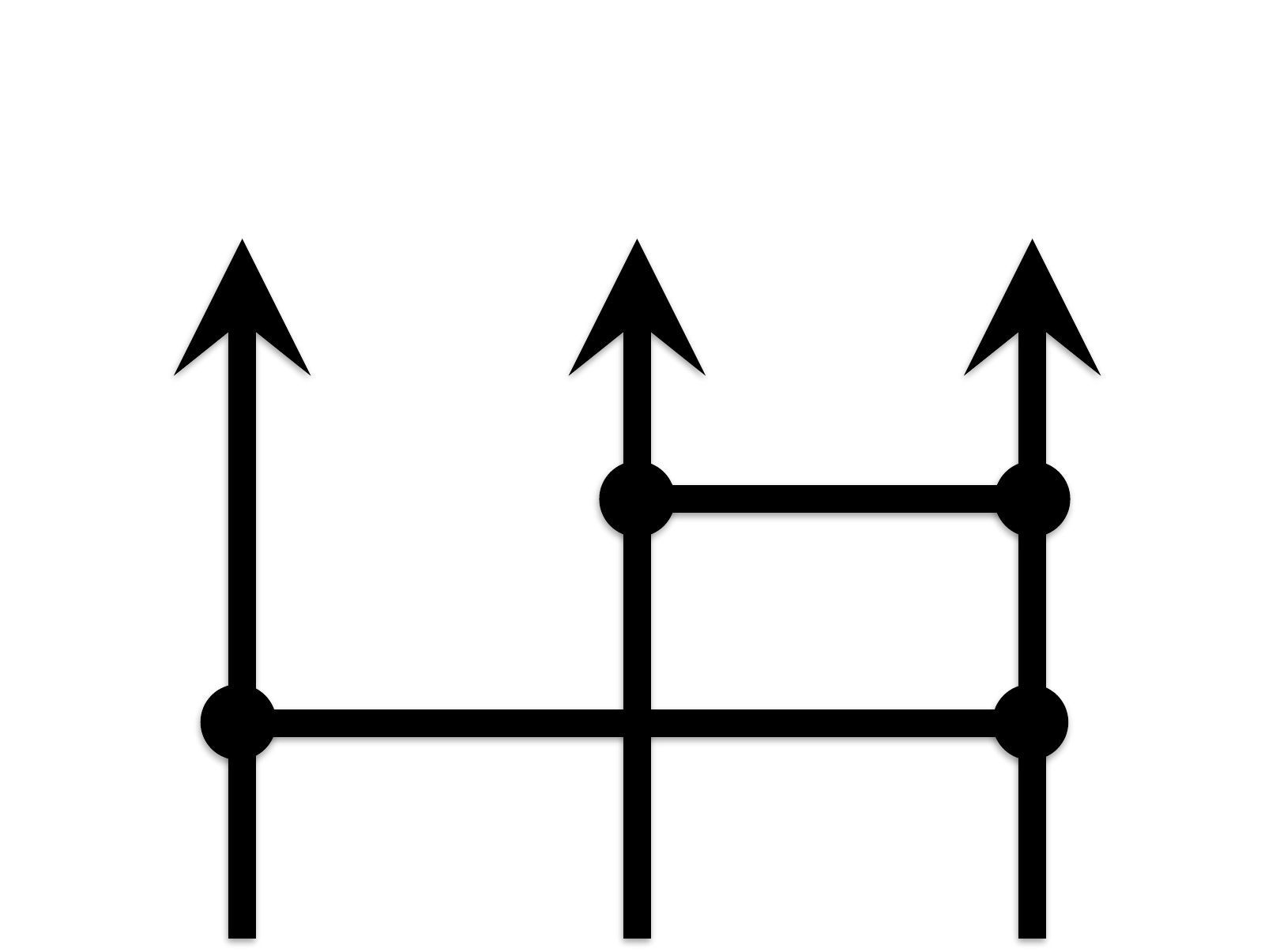}}}
-
\hspace*{-3pt}\raisebox{-.3\height}{\scalebox{0.1}{\includegraphics{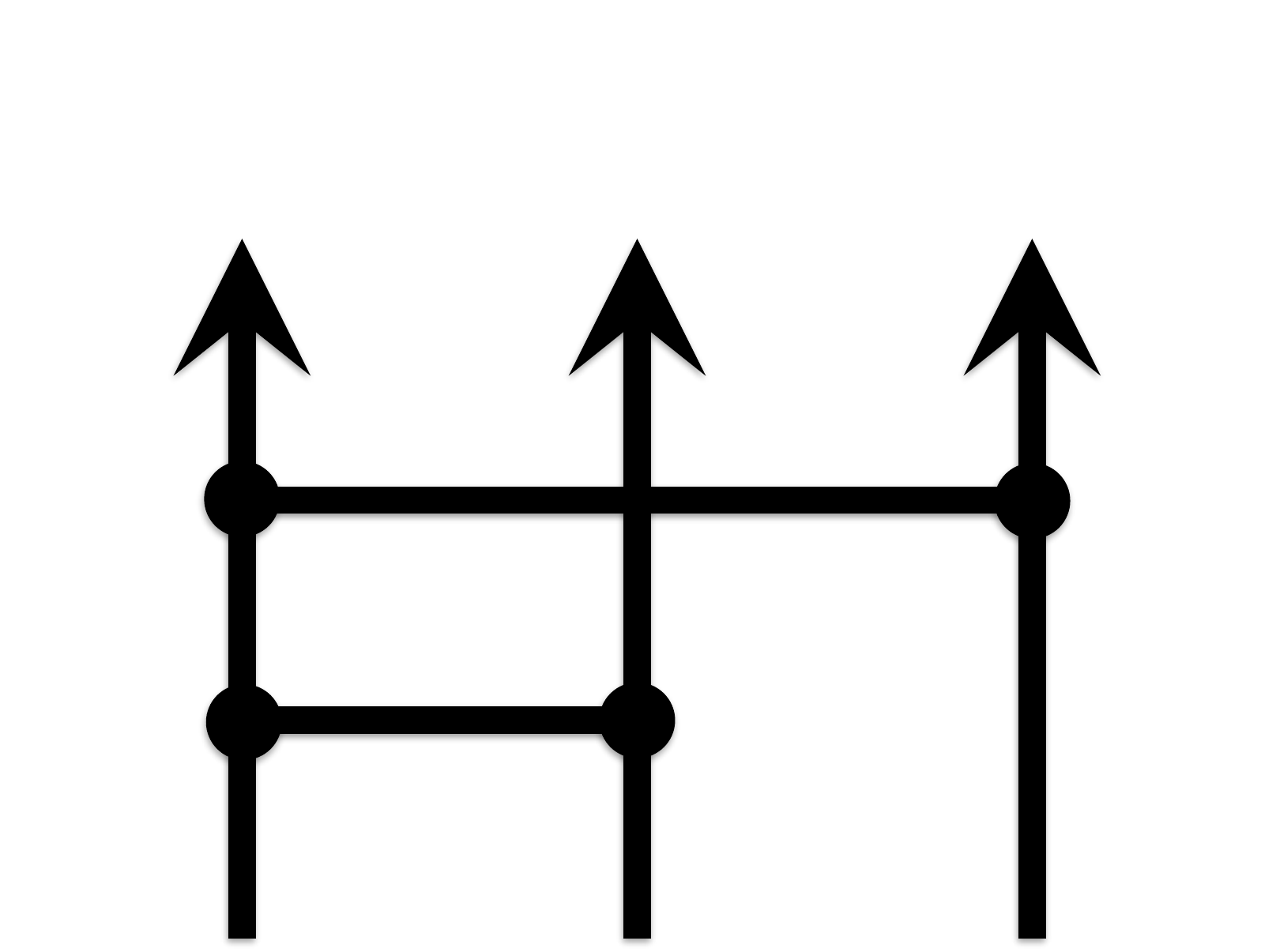}}}
-
\hspace*{-3pt}\raisebox{-.3\height}{\scalebox{0.1}{\includegraphics{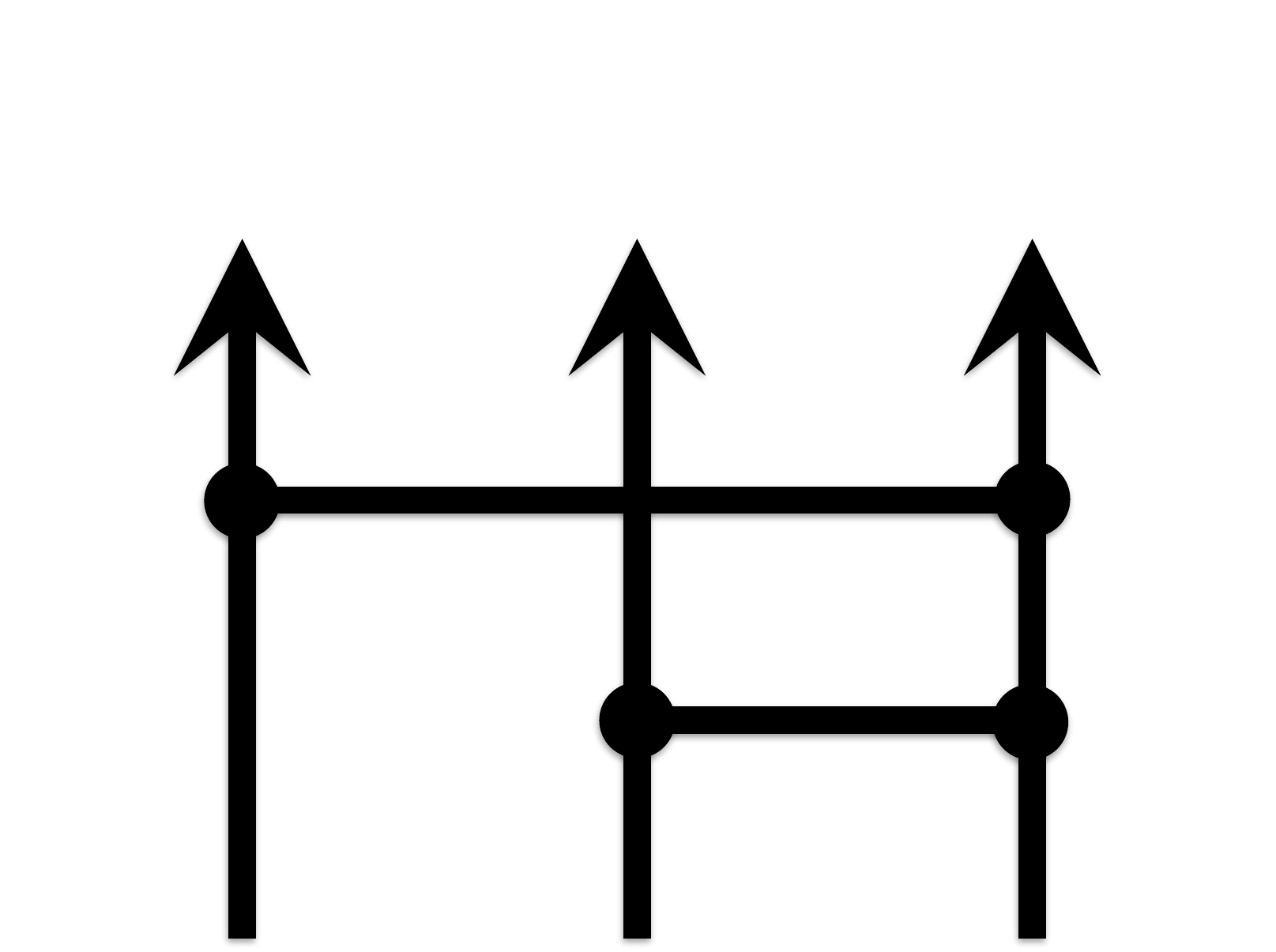}}}
\}.
$$
(Note the difference between a vertex, indicated by a dot, and a crossing of edges as an effect
of the planarity of the drawing.)
Then the functions $f:\GG\to\oF$ annihilating $Q$ are the `weight systems'.
Moreover, $\hat p_R(Q)=0$ if and only if $R$ comes from a representation of a Lie algebra
(cf.\ [3]).

\smallskip
\noindent
{\em Groups.}
Define for a group $\Gamma$:
$$
Q:=
\{
\hspace*{-12pt}\raisebox{-.33\height}{\scalebox{0.12}{\includegraphics{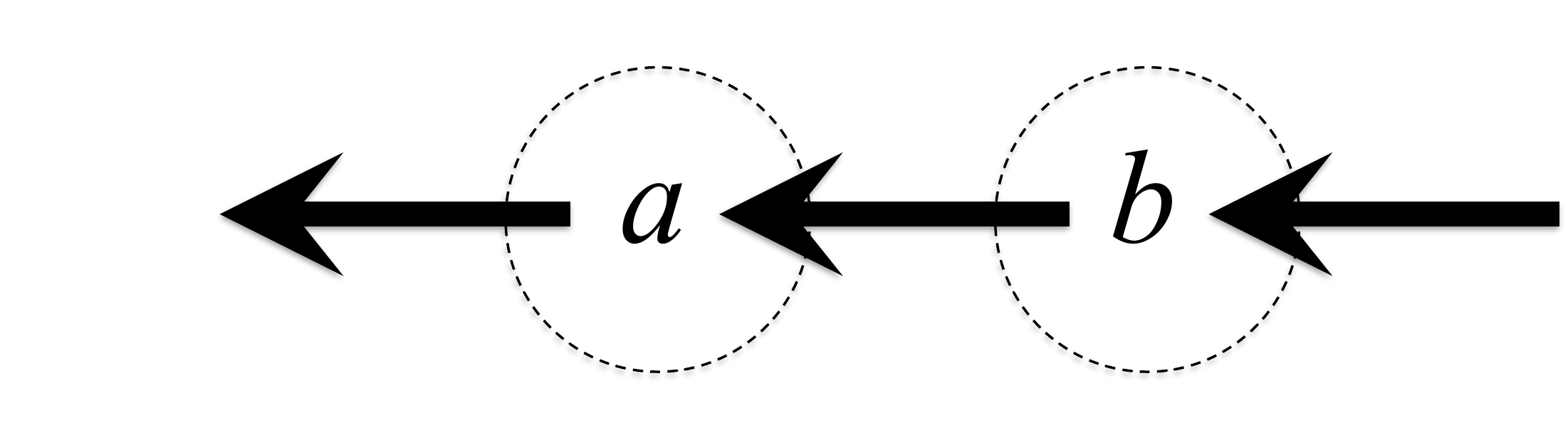}}}
-
\hspace*{-11pt}\raisebox{-.33\height}{\scalebox{0.12}{\includegraphics{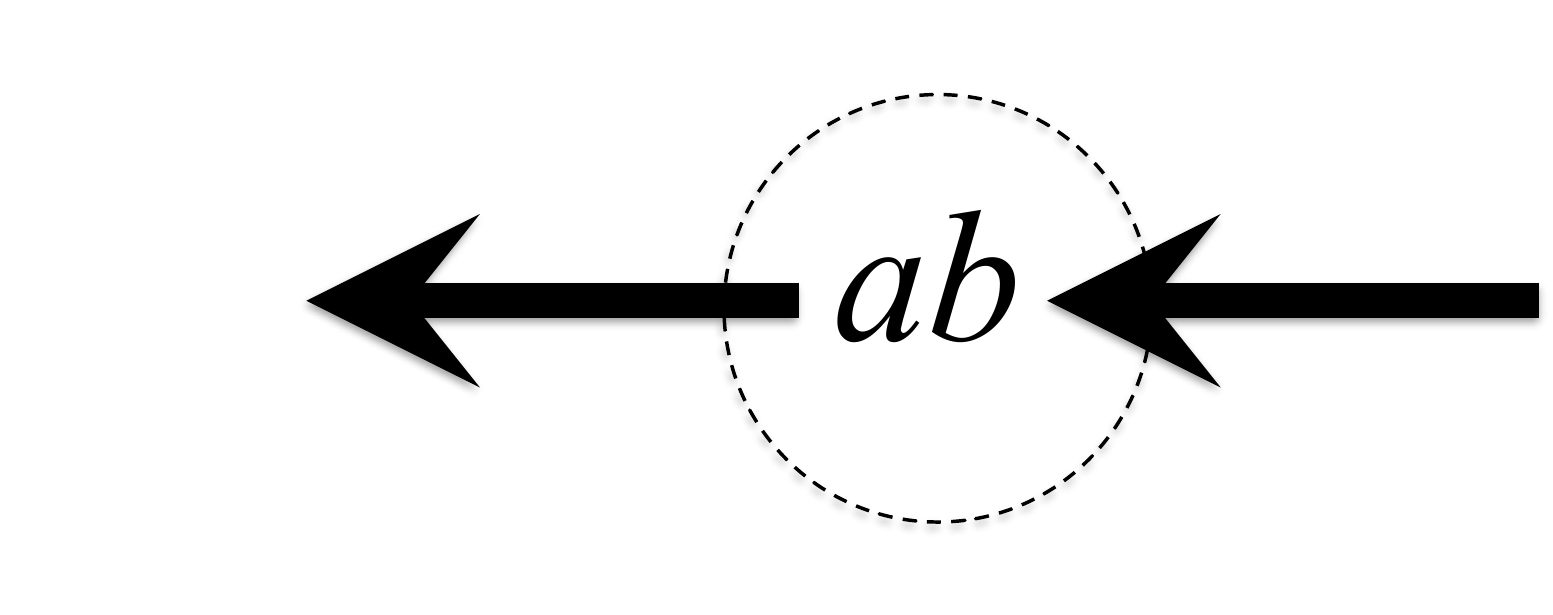}}}
\mid
a,b\in \Gamma\}\cup
\{
\hspace*{-10pt}\raisebox{-.33\height}{\scalebox{0.12}{\includegraphics{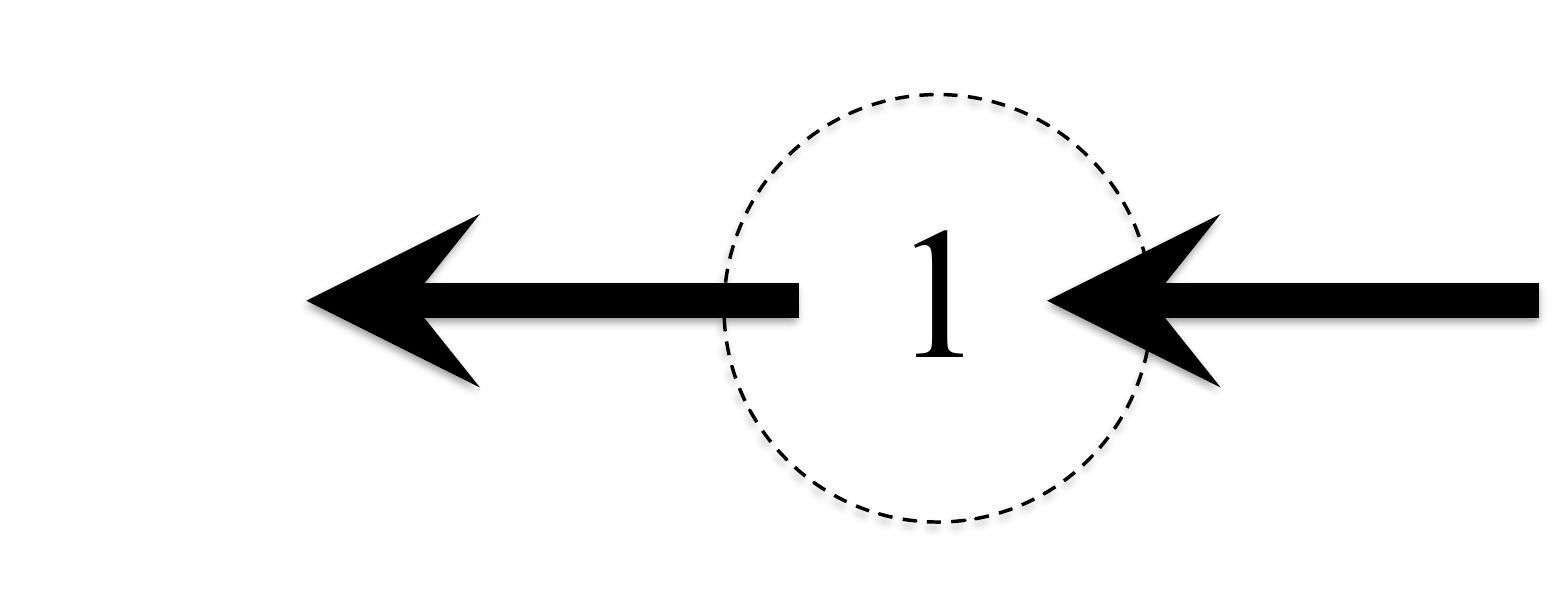}}}
-
\hspace*{-9pt}\raisebox{-.3\height}{\scalebox{0.12}{\includegraphics{td_id.pdf}}}
\}
,
$$
where $1$ is the unit of $\Gamma$.
Then $\hat p_R(Q)=0$ if and only if $R$ is a representation of $\Gamma$.

\smallskip
\noindent
{\em Algebra template.}
Define
$$Q:=
\{
\hspace*{-9pt}\raisebox{-.3\height}{\scalebox{0.09}{\includegraphics{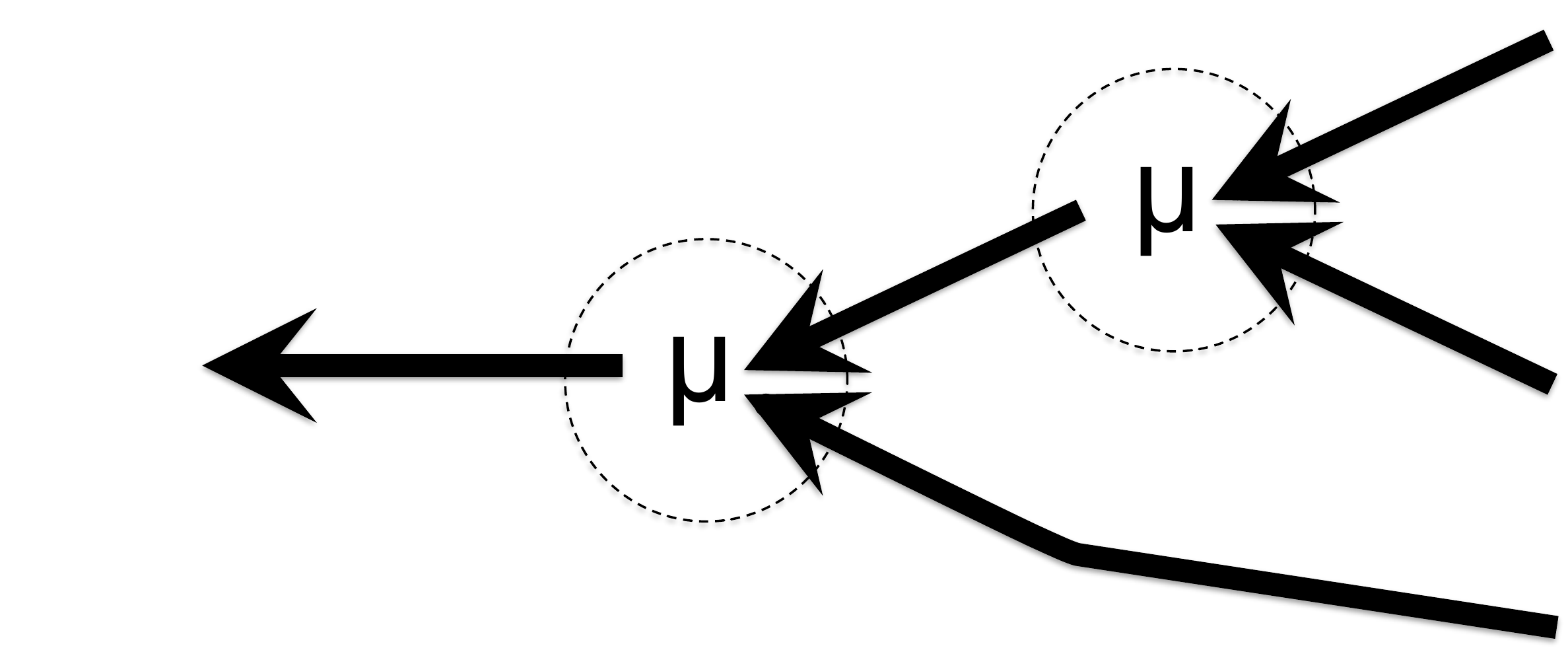}}}
-
\hspace*{-9pt}\raisebox{-.3\height}{\scalebox{0.09}{\includegraphics{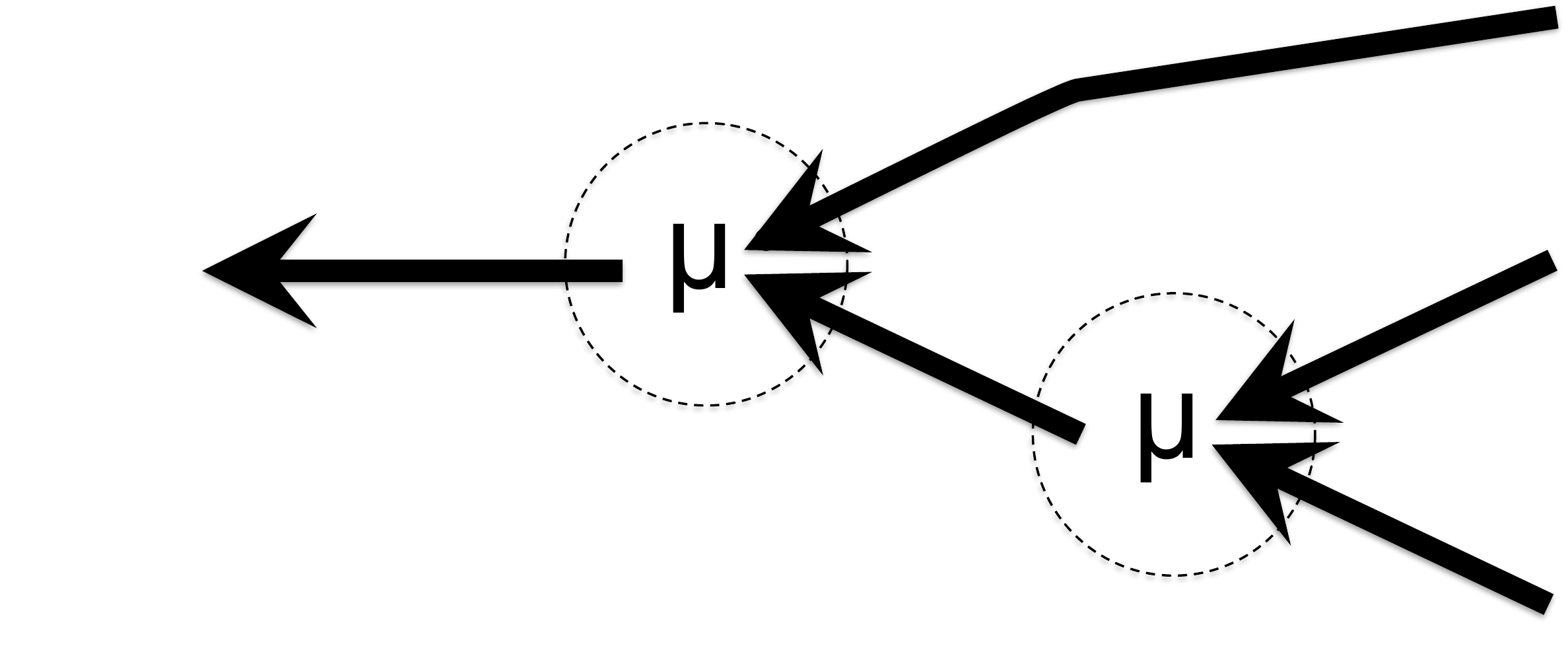}}}
,
\hspace*{-9pt}\raisebox{-.2\height}{\scalebox{0.09}{\includegraphics{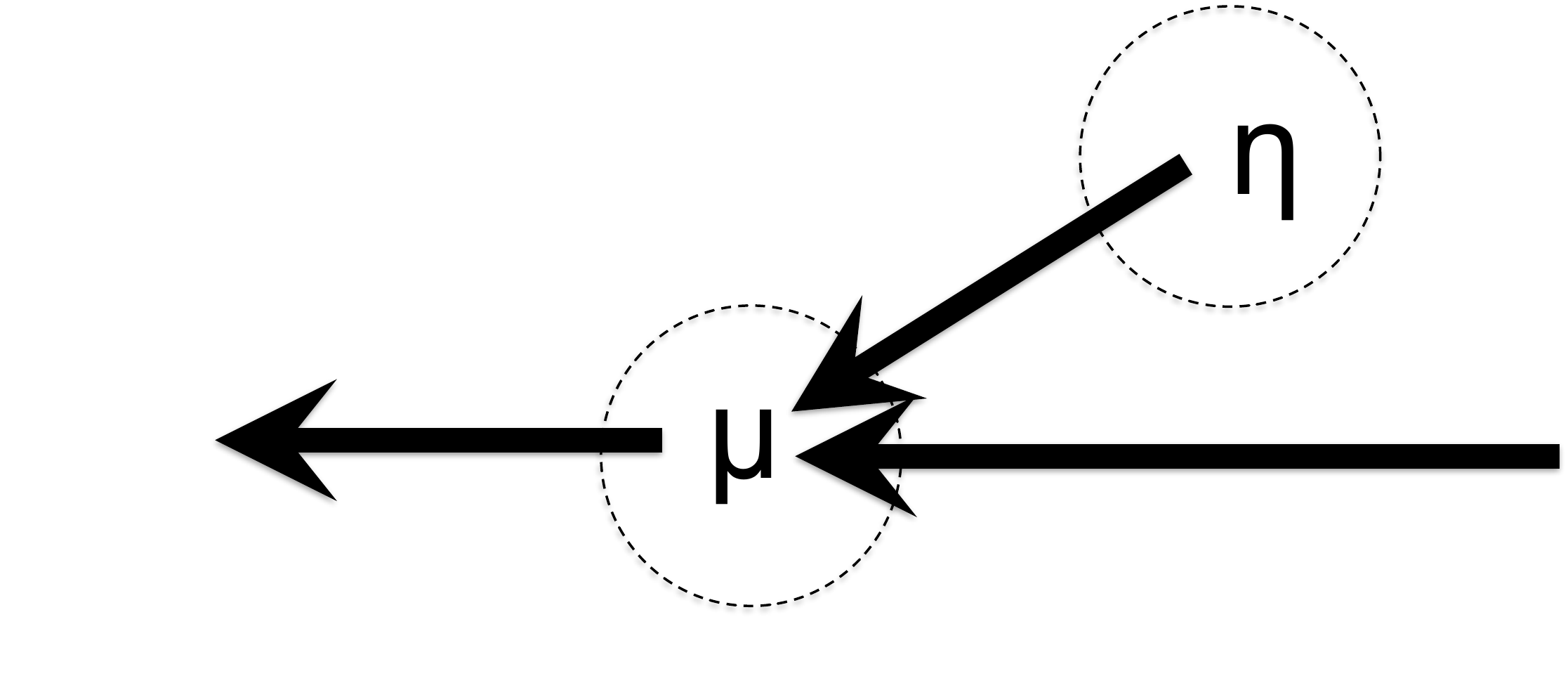}}}
-
\hspace*{-9pt}\raisebox{-.3\height}{\scalebox{0.09}{\includegraphics{td_id.pdf}}}
,
\hspace*{-9pt}\raisebox{-.55\height}{\scalebox{0.09}{\includegraphics{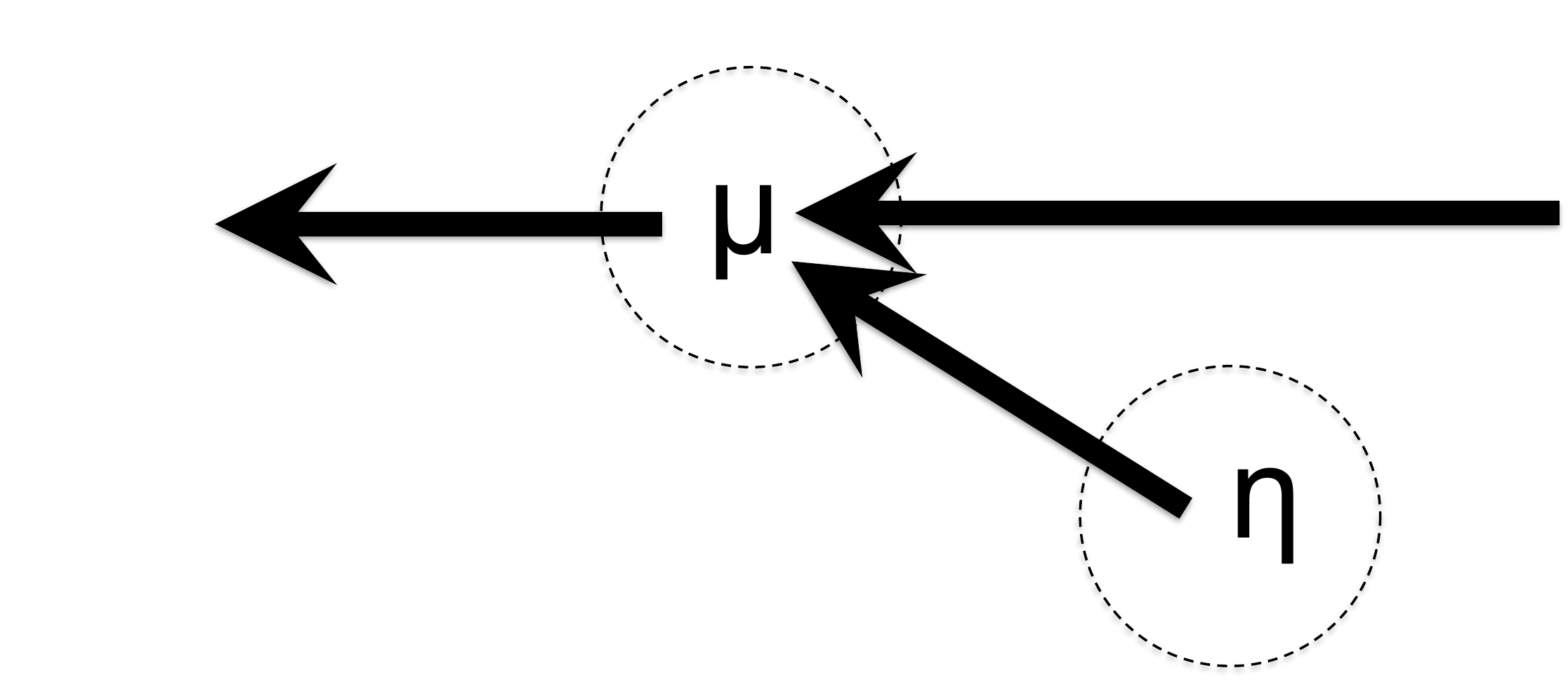}}}
-
\hspace*{-9pt}\raisebox{-.3\height}{\scalebox{0.09}{\includegraphics{td_id.pdf}}}
\}.$$
Then $\hat p_R(Q)=0$ if and only if $R(\mu)$ and $R(\eta)$ are the multiplication tensor and the unit
of a finite-dimensional unital associative $\oF$-algebra.

\smallskip
\noindent
{\em Hopf algebra template.}
The Hopf algebra axioms can similarly be translated into quantum diagrams (cf.\ [12]).

\smallskip
\noindent
{\em Directed graphs.}
Let $Q$ be the collection of all quantum webs
$$
\hspace*{-9pt}\raisebox{-.4\height}{\scalebox{0.15}{\includegraphics{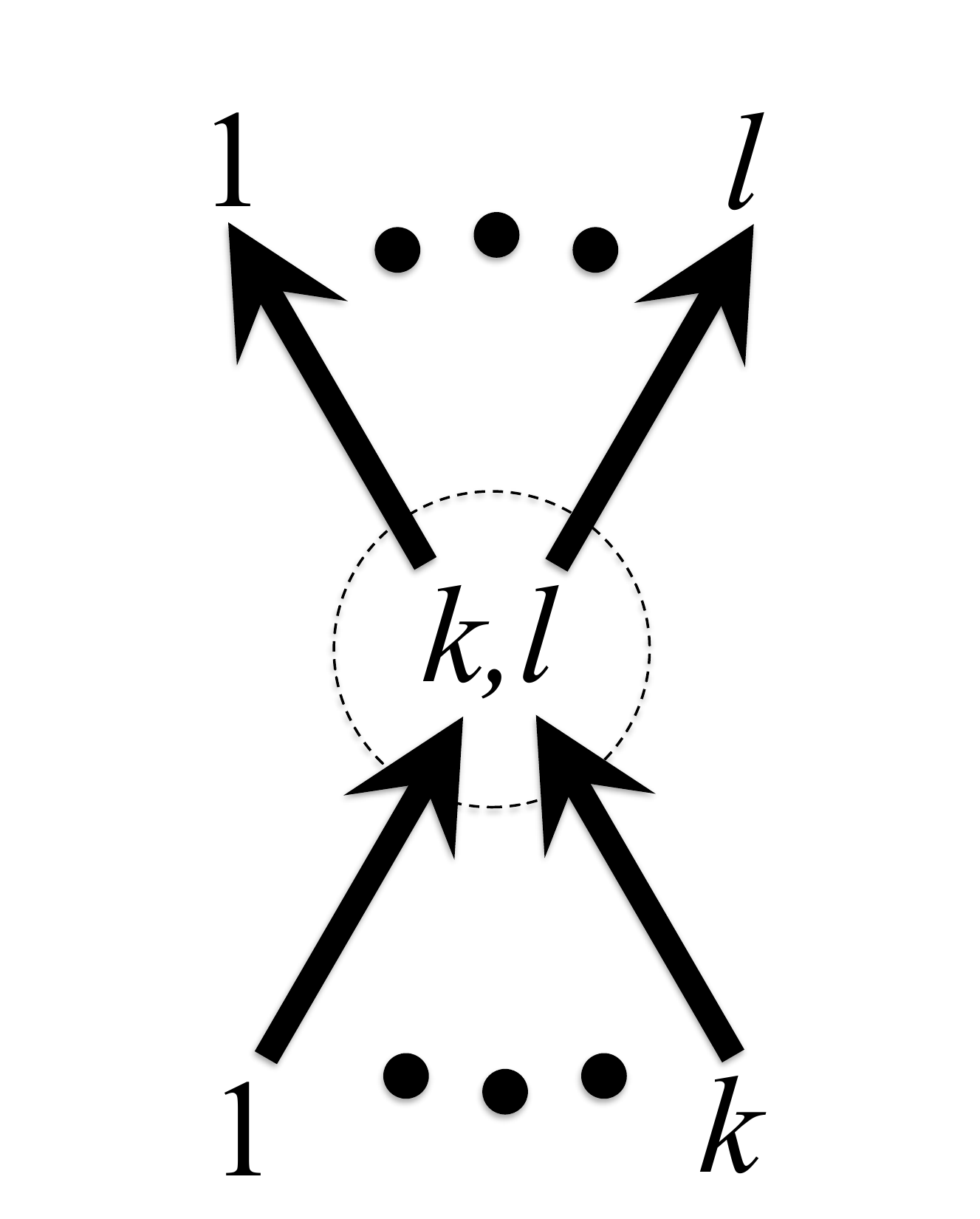}}}
-
\hspace*{-9pt}\raisebox{-.4\height}{\scalebox{0.15}{\includegraphics{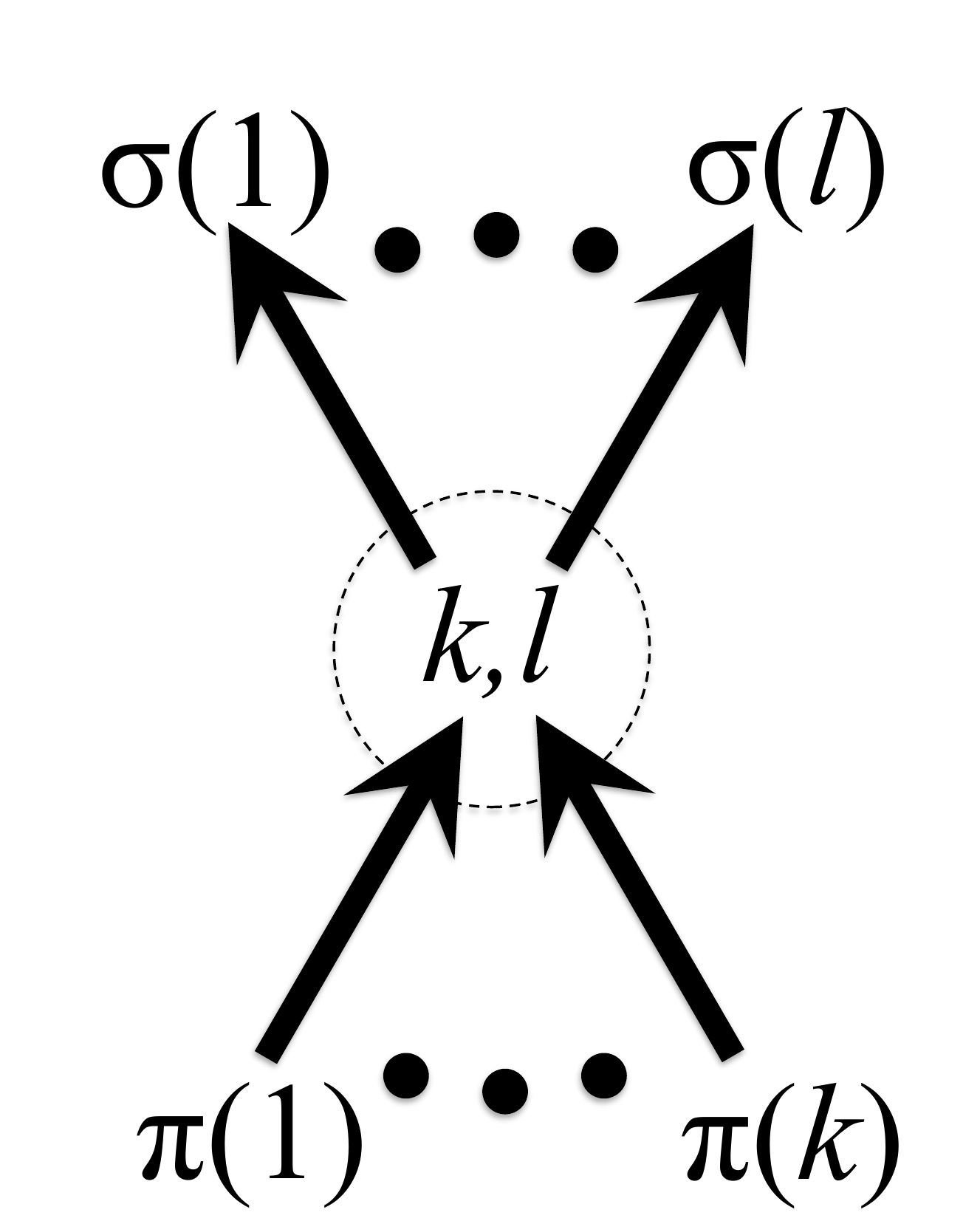}}}
$$
with $k,l\in\oZ_+$ and $\pi\in S_k$, $\sigma\in S_l$.
Then $\hat p_R(Q)=0$ amounts to requiring that $R$ is symmetric under permutations of entering
edges and under permutations of leaving edges.
Thus we deal with invariants of ordinary directed graphs, with no ordering of edges.
This case was considered in [4], and Theorem \ref{13me14a} forms a generalization of its
result.

\tussenkop{Nondegeneracy}

In these examples we were considering tensor representations $R$ with $p_R=f$ and $\hat p_R(Q)=0$, where $Q$
is some given collection of quantum webs annihilated by $f$.
In fact, for each trace $f$ there exists a tensor representation $R$ with $f=p_R$ such that
$\hat p_R(\omega)=0$ for {\em each} quantum web $\omega$ annihilated by $f$.
To describe this more precisely, we define nondegeneracy of tensor representations.

Call a tensor representation $R:T\to \mbox{\rm T}(V)$ {\em nondegenerate} if for each
quantum web $\omega$ with $\hat p_R(\omega)\neq 0$, there exists
$W\in\WW$ with $p_R(\omega\cdot W)\neq 0$.
In other words, in view of \rf{5ja15a}, $R$ being nondegenerate means that the subspace $\hat p_R(\oF\WW)$ of $\mbox{\rm T}(V)$
is nondegenerate with respect to the standard bilinear form on $\mbox{\rm T}(V)$.
Or: $\hat p_R(\omega)=0$ for each quantum web $\omega$ annihilated by $p_R$.

Call $R$ {\em strongly nondegenerate} if each finite $U\subseteq T$ is contained in some finite
$S\subseteq T$ with the restriction $R|S$ of $R$ to $S$ being nondegenerate.
The proof of the theorem below implies that this is equivalent to:
there is a finite subset $U\subseteq T$ such that $R|S$ is nondegenerate for each finite $S$
with $U\subseteq S\subseteq T$.
So strong nondegeneracy implies nondegeneracy, and if $T$ is finite, the two concepts coincide.
(Actually, we have no example of a nondegenerate $R$ which is not strongly nondegenerate.)

In the following theorem, `unique' means: up to the natural action of $\GL(V)$ on the set of tensor
representations $T\to \mbox{\rm T}(V)$.

\thm{27de14b}{
For each trace $f$ there exists a unique strongly nondegenerate tensor representation $R$ with
$p_R=f$.
}

\bigskip
For given sets $T$ of types and $Q$ of quantum webs,
it is a fundamental question to determine the collection $\overline Q$ of quantum webs $\omega$
that are annihilated by each function $f:\GG_T\to\oF$ annihilating $Q$.
For the virtual link diagram example this contains the question which virtual link diagrams are
equivalent under Reidemeister moves.

A related question is whether for each $f$ annihilating $Q$ and each quantum diagram $\gamma$
with $f(\gamma)\neq 0$, there exists a trace $p_R$ annihilating $Q$ with $p_R(\gamma)\neq 0$ (`detecting $\gamma$').
For instance, for the multiloop chord diagram example, this question was answered negatively by
Vogel [22].

\sectz{Some applications of invariant theory}

We give a few consequences of invariant theory, as preparation to the proof of Theorems
\ref{13me14a} and \ref{27de14b} in Section \ref{21no10a}.
In this section, fix a finite-dimensional linear space $V$.

If $T$ is finite, then $\RR_T$ is a finite-dimensional linear space, which can be described as:
\dyy{ster}{
\RR_T=\bigoplus_{t\in T}V^{*\otimes\iota(t)}\otimes V^{\otimes o(t)}.
}
Then the following is a direct application of the first fundamental theorem (FFT) of invariant
theory for $\GL(V)$ (cf.\ [6] Corollary 5.3.2), where as usual $\OO(Z)$
denotes the set of regular $\oF$-valued functions on a variety $Z$, while if $\GL(V)$ acts on
a set $S$, then $S^{\GL(V)}$ is the set of $\GL(V)$-invariant elements of $S$:
\dyy{30de14a}{
\hat p(\oF\WW)=(\OO(\RR_T)\otimes \mbox{\rm T}(V))^{\GL(V)}.
}

\vspace{4mm}
\noindent
{\bf Proposition \thebewering.}\label{22no08a}{\it
Let $T$ be finite and $R\in\RR_T$.
Then $R$ is nondegenerate if and only if the orbit $\GL(V)\cdot R$ is closed.
}

\pf
{\em Sufficiency.}
Let $C:=\GL(V)\cdot R$ be closed and let $\omega\in\oF\WW$ be such that $p_R(\omega\cdot W)=0$ for each
$W\in\WW$.
Suppose that $\hat p_R(\omega)\neq 0$.
As the function $\hat p(\omega)$ is $\GL(V)$-equivariant, $\hat p_{R'}(\omega)\neq 0$ for all $R'\in C$.
Hence, since $C$ is closed, by the Nullstellensatz there exists $q\in\OO(\RR_T)\otimes \mbox{\rm T}(V)$ with
$\hat p(\omega)(R')\cdot q(R')=1$ for each $R'\in C$.
Applying the Reynolds operator, we can assume that $q$ is $\GL(V)$-equivariant
(as $\hat p(\omega)$ is $\GL(V)$-equivariant).
So by \rf{30de14a}, $q=\hat p(\omega')$ for some quantum web $\omega'$.
Then $1=\hat p(\omega)(R)\cdot\hat p(\omega')(R)=p(\omega\cdot\omega')(R)=p_R(\omega\cdot\omega')=0$,
a contradiction.

{\em Necessity.}
Let $F:=\{R'\in\RR_T\mid p_{R'}=p_R\}$.
So $F$ is the set of all $R'\in\RR_T$ with $d(R')=d(R)$ for each $\GL(V)$-invariant regular
function $d$ on $\RR_T$ (by \rf{30de14a}).
Hence $F$ is a fiber of the projection $\RR_T\to\RR_T//\GL(V)$.
So $F$ contains a unique closed $\GL(V)$-orbit $C$ ([2]).

Suppose $R\not\in C$.
Then there exists $q\in\OO(\RR_T)$ with $q(C)=0$ and $q(R)\neq 0$.
Let $U$ be the $\GL(V)$-module spanned by $\GL(V)\cdot q$.
The morphism $\phi:\RR_T\to U^*$ with $\phi(R')(u)=u(R')$ (for
$R'\in\RR_T$ and $u\in U$) is $\GL(V)$-equivariant.

Let $\varepsilon:U^*\to \mbox{\rm T}(V)$ be an embedding of $U^*$ as $\GL(V)$-submodule of $\mbox{\rm T}(V)$.
(This exists, as $U$ is spanned by a $\GL(V)$-orbit, so that each irreducible $\GL(V)$-module occurs with
multiplicity at most 1 in $U$.)
So $\varepsilon\circ\phi$ is a $\GL(V)$-equivariant morphism $\RR_T\to \mbox{\rm T}(V)$.
In other words, $\varepsilon\circ\phi$ belongs to $(\OO(\RR_T)\otimes \mbox{\rm T}(V))^{\GL(V)}$,
which is by \rf{30de14a} equal to $\hat p(\oF\WW)$.
Hence $\varepsilon\circ\phi=\hat p(\omega)$ for some $\omega\in\oF\WW$.
As $\phi(R)\neq 0$ (since $\phi(R)(q)=q(R)\neq 0$), we have $\hat p(\omega)(R)\neq 0$.
As $R$ is nondegenerate, there is a web $W$ with $p(\omega\cdot W)(R)\neq 0$.
So $p_R(\omega\cdot W)\neq 0$.
However, for any $R'\in C$, $p(\omega\cdot W)(R')=\hat p(\omega)(R')\cdot\hat p(W)(R')=0$,
since $\hat p(\omega)(R')=\varepsilon\circ\phi(R')=0$, as $q(C)=0$.
So $p_{R'}(\omega\cdot W)=0$ while
$p_{R}(\omega\cdot W)\neq 0$, contradicting the fact that
$p_{R'}(G)=p_{R}(G)$ for each diagram $G$.
\bx

\sect{21no10a}{Proof of Theorems \ref{13me14a} and \ref{27de14b}}

I.
We first show necessity in Theorem \ref{13me14a}.
Let $f=p_R$ for some tensor representation $R:T\to \mbox{\rm T}(V)$, where $V$ is a $d$-dimensional linear space
with $d\leq n$.
Clearly, $p_R$ is multiplicative.
Moreover, $\hat p_R(\Delta_{n+1})=0$.
Indeed, consider an `edge coloring' $\phi$ in the summation \rf{5ja15b}.
As $d<n+1$, two edges of $\Delta_{n+1}$ have the same $\phi$-value, say edges $e_i$ and $e_j$ ($i\neq j$).
Let $\sigma$ be the permutation in $S_{n+1}$ swapping $i$ and $j$.
Then $J_{\pi}$ and $J_{\sigma\circ\pi}$ cancel each other out for this $\phi$, as
$\pi$ and $\sigma\circ\pi$ have opposite signs.
So for each fixed $\phi$, the term in \rf{5ja15b} is 0.
Therefore, $\hat p_R(\Delta_{n+1})=0$, hence $p_R(\Delta_{n+1}\cdot W)=0$ for each
$W\in\WW$, by \rf{5ja15a}.

\medskip
\noindent
II.
We next show that the condition in Theorem \ref{13me14a} implies the existence of
a strongly nondegenerate tensor representation $R$ with $p_R=f$.

Let $f:\GG\to\oF$ be multiplicative and annihilate $\Delta_{n+1}$.
We can assume that $n$ is smallest with this property.
Then:
$$
f(\loop)=n.
$$
Indeed, by the minimality of $n$, there exists $W\in\WW_{n,n}$ with $f(\Delta_n\cdot W)\neq 0$.
Let $W'\in\WW_{n+1,n+1}$ be obtained from $W$ by adding one directed edge disjoint from $W$,
with both ends labeled $n+1$.
Then $0=f(\Delta_{n+1}\cdot W')=(f(\loop)-n)f(\Delta_n\cdot W)$.
So $f(\loop)=n$.

From now on in this proof, fix an $n$-dimensional $\oF$-linear space $V$.
So $\RR_T$ and $p:\oF\GG\to\OO(\RR_T)$ are well-defined.
Then $p$ is an algebra homomorphism, with respect to the $\cdot$ product on the space $\oF\GG$
of quantum diagrams (which is for diagrams just the disjoint union).

\bigskip
\noindent
{\em Claim.}
 {\it
$\Ker~p\subseteq\Delta_{n+1}\cdot\oF\WW$.
}

\pfcl
Let $\gamma\in\oF\GG$ with $p(\gamma)=0$.
We prove that $\gamma\in\Delta_{n+1}\cdot\oF\WW$.
By splitting $p(\gamma)$ into homogeneous components,
we can assume that $\gamma$ is a linear combination of diagrams that all have the same number $d$
of vertices; and more strongly, that
all have the same number $m_t$ of vertices of type $t$, for any $t\in T$.
Then we can assume (by renaming and deleting unused types) that $T:=\{1,\ldots,a\}$ for some $a$,
$m_1+\cdots+m_a=d$, and $m_t\geq 1$ for all $t\in T$.

In fact, we can assume that $m_t=1$ for each $t\in T$.
To see this, consider a type $t\in T$ with $m_t\geq 2$, and introduce a new type, named $a+1$,
with $\iota(a+1)=\iota(t)$ and $o(a+1)=o(t)$.
For any diagram $G$, let $G'$ be the sum of those diagrams that can be obtained
from $G$ by changing the type of one vertex of type $t$ to type $a+1$.
So $G'$ is the sum of $m_t$ diagrams.

To describe $p(G')$, let $B$ be a basis of term $V^{*\otimes \iota(t)}\otimes V^{\otimes o(t)}$ in \rf{ster}.
For $b\in B$, let $b^*$ be the corresponding element in the basis dual to $B$,
and let $b'$ be the element corresponding to $b^*$
for the new term $V^{*\otimes\iota(a+1)}\otimes V^{\otimes o(a+1)}$ in \rf{ster}.
Then
\dyy{3ja15a}{
p(G')=\sum_{b\in B}\frac{d}{db^*}p(G)b'.
}
Let $\gamma'$ be obtained by replacing each $G$ in $\gamma$ by $G'$.
As $p(\gamma)=0$, \rf{3ja15a} gives $p(\gamma')=0$.
Moreover, $\gamma'\in\Delta_{n+1}\cdot\oF\WW$ implies
$\gamma\in\Delta_{n+1}\cdot\oF\WW$, as we can apply a reverse map $G'\mapsto G$
(where $G$ is obtained from $G'$ by replacing type $a+1$ by $t$ and dividing by $m_t$).

Repeating this operation we finally obtain that each type occurs precisely once in each diagram occurring in $\gamma$.
So finally $a=d$.
Then $\sum_{t=1}^a\iota(t)=\sum_{t=1}^ao(t)=:m$.

Make the following web $W\in\WW_{m,m}$, having vertices $v_1,\ldots,v_a$, where $v_t$ has type $t$, for $t\in T$,
and having in addition $m$ roots and $m$ sinks.
The tails of the edges entering $v_t$ are roots, labeled $\overline\iota(t)+1,\ldots,\overline\iota(t)+\iota(t)$, in order,
where $\overline\iota(t):=\sum_{i<t}\iota(i)$.
Moreover, the heads of the edges leaving $v_t$ are sinks, labeled $\overline o(t)+1,\ldots,\overline o(t)+o(t)$, in order,
where $\overline o(t):=\sum_{i<t}o(i)$.

For each $\pi\in S_m$, let $G_{\pi}:=J_{\pi}\cdot W$.
Then for each diagram $G$ with $a$ vertices, of types $1,\ldots,a$ respectively,
there exists a unique $\pi\in S_m$ with $G=G_{\pi}$.
Moreover, for all $y_1,\ldots,y_m\in V^*$ and $z_1,\ldots,z_m\in V$:
\dyy{19ja15a}{
\prod_{i=1}^my_{\pi(i)}(z_i)=p(G_{\pi})\big(
\bigoplus_{t=1}^a
\big(
\bigotimes_{i=1}^{\iota(t)}y_{\overline\iota(t)+i}\otimes\bigotimes_{j=1}^{o(t)}z_{\overline o(t)+j}
\big)
\big).
}

We can write for unique $\lambda_{\pi}\in\oF$ (for $\pi\in S_m)$:
$$
\gamma=\sum_{\pi\in S_m}\lambda_{\pi}G_{\pi}.
$$
Define the following polynomial $q\in\OO(\oF^{m\times m})$:
$$
q(X):=\sum_{\pi\in S_m}\lambda_{\pi}\prod_{i=1}^{m}x_{\pi(i),i}
$$
for $X=(x_{i,j})_{i,j=1}^m\in\oF^{m\times m}$.
Note that $q$ determines $\gamma$.
As $p(\gamma)=0$, \rf{19ja15a} implies that $q((y_i(z_j))_{i,j=1}^{m})$
for all $y_1,\ldots,y_m\in V^*$ and $z_1,\ldots,z_m\in V$.
This implies, by the second fundamental theorem (SFT) of invariant theory for $\GL(V)$
(cf.\ [6] Theorem 12.2.12),
that $q$ belongs to the ideal in $\OO(\oF^{m\times m})$ generated by the
$(n+1)\times(n+1)$ minors of $\oF^{m\times m}$.
That is,
$$
q=\sum_{I,J\subseteq [m]\atop|I|=|J|=n+1}q_{I,J}\det(X_{I,J}),
$$
where $X_{I,J}$ is the $I\times J$ submatrix of $X\in\oF^{m\times m}$ and $q_{I,J}$ belongs
to $\OO(\oF^{m\times m})$.
As each term of $q$ comes from a permutation, the variables in any row of $X$ have total degree
1 in $q$.
Similarly, the variables in any column of $X$ have total degree 1 in $q$.
This implies that we can assume that each $q_{I,J}$ has total degree 1 in rows of $X$ with index not in $I$
and total degree 0 in rows in $X$ with index in $I$.
Similarly for columns with respect to $J$.
This implies $\gamma\in\Delta_{n+1}\cdot\oF\WW$.
\openbx

This claim and the condition in Theorem \ref{13me14a} imply that $\Ker~p\subseteq\Ker f$.
Hence there exists a linear function $\hat f:p(\oF\GG)\to\oF$ such that $\hat f\circ p=f$.
Then $\hat f$ is a unital algebra homomorphism,
as $\hat f(1)=\hat f(p(\emptyset))=f(\emptyset)=1$, and as for all $G,H\in\GG$:
$$
\hat f(p(G)p(H))=\hat f(p(GH))=f(GH)=f(G)f(H)=\hat f(p(G))\hat f(p(H)).
$$

Now first suppose that $T$ is finite.
Since $p(\oF\GG)=\OO(\RR_T)^{\GL(V)}$ by \rf{ster},
there exists $R\in\RR_T$ with $\hat f(q)=q(R)$ for all $q\in p(\oF\GG)$ (by the Nullstellensatz).
So
$$
p_R(G)=p(G)(R)=\hat f(p(G))=f(G),
$$
for all $G\in\GG_T$,
proving Theorem \ref{13me14a} for finite $T$.
By the closed orbit theorem ([2]), we can assume that the orbit $\GL(V)\cdot R$ is closed.
Then, by Proposition \ref{22no08a}, $R$ is nondegenerate.

Suppose next that $T$ is infinite.
Consider any finite subset $U$ of $T$.
Let $F_U$ be the variety of tensor representations $R:U\to \mbox{\rm T}(V)$ such that $p_R=f|\GG_U$.
We saw above that, as $U$ is finite, $F_U\neq\emptyset$.
In fact, $F_U$ is a fiber of the projection $\RR_U\to \RR_U//\GL(V)$.
Hence $F_U$ contains a unique $\GL(V)$-orbit $C_U$ of minimal (Krull) dimension
(cf.\ [2] 1.11 and 1.24).
(It is in fact the unique closed orbit in $F_U$.)

Since $\dim(C_U)\leq\dim(\GL(V))$ for each finite $U\subseteq T$, there exists a finite $U\subseteq T$ with $\dim(C_U)$ as large as possible.
Then for each finite $S\subseteq T$ with $U\subseteq S$:
\dy{29de14a}{
$\dim(C_S)=\dim(C_U)$ and $\pi_{S,U}(C_S)=C_U$,
}
where $\pi_{S,U}$ is the natural projection $\RR_S\to \RR_U$.
Indeed,
\dyy{15me14b}{
\dim(C_U)\leq\dim(\pi_{S,U}(C_S))\leq\dim(C_S).
}
To see this, note that $\pi_{S,U}(C_S)\subseteq F_U$.
Then the first inequality in \rf{15me14b} follows from the fact that $\pi_{S,U}(C_S)$ is a
$\GL(V)$-orbit
in $F_U$, and that $C_U$ has minimal dimension among all $\GL(V)$-orbits in $F_U$.

By the maximality of $\dim(C_U)$, we have equality throughout in \rf{15me14b}.
As $\pi_{S,U}(C_S)$ is a $\GL(V)$-orbit and as $C_U$ is the unique orbit in $F_U$ of minimal dimension, \rf{29de14a} follows.

Choose an arbitrary $R\in C_U$ and consider some finite $S\supseteq U$.
We extend $U$ and $R$ if possible as follows.
By \rf{29de14a}, there exists at least one $R'\in C_S$ with $\pi_{S,U}(R')=R$.
As $\pi_{S,U}$ is $\GL(V)$-equivariant, for the stabilizers one has
\dyy{15me14d}{
\GL(V)_{R'}\subseteq\GL(V)_R.
}
Suppose that there exists $R''\neq R'$ in $C_S$ with $\pi_{S,U}(R'')=R$.
So some $g\in\GL(V)$ moves $R'$ to $R''$, while it leaves $R$ invariant.
Then we have strict inclusion in \rf{15me14d}.
Now replace $U,R$ by $S,R'$.

In respect of \rf{15me14d}, the finite basis theorem implies that we can do
such replacements only a finite number of times.
So we end up with a finite $U$ and $R\in C_U$ such that $\pi_{S,U}$ is injective on $C_S$.
Hence for each $t\in T$ there is a unique $R_t\in C_{U\cup\{t\}}$
such that $R_t|U=R$.

Define a tensor representation $P:T\to\mbox{\rm T}(V)$ by $P(t):=R_t(t)$ for $t\in T$.
This implies that for each finite $U'\supseteq U$, $P|U'$ belongs to $C_{U'}$.
Hence by Proposition \ref{22no08a}, $P|U'$ is nondegenerate.
Concluding, $P$ is strongly nondegenerate.

\medskip
\noindent
III.
Finally we show the uniqueness of a strongly nondegenerate tensor representation $R$ with $p_R=f$, up to the
action of $\GL(V)$ on $\RR_T$.
Let $R$ and $R'$ be strongly nondegenerate tensor representations with $f=p_R=p_{R'}$.
It suffices to show that for each finite $S\subseteq T$ there exists $g\in\GL(V)$ such that
$R'|S=g\cdot R|S$ (since then, by the finite basis theorem,
we can choose a finite $S\subseteq T$ with the variety
$\Gamma_S:=\{g\in\GL(V)\mid R'|S=g\cdot R|S\}$ minimal among all $\Gamma_{S'}$ with finite $S'\supseteq S$,
implying that for each $t\in T$ and $g\in\Gamma_S$ one has $R'|S\cup\{t\}=g\cdot R|S\cup\{t\}$).

Let again $U\subseteq T$ be finite with $\dim(C_U)$ maximal.
Hence it suffices to show that for each finite $U'\subseteq T$ with $U'\supseteq U$,
$\GL(V)\cdot R|U'=\GL(V)\cdot R'|U'$.

As $U\subseteq U'$, $\dim(C_{U'})$ is maximal (by \rf{29de14a}).
As $R$ and $R'$ are strongly nondegenerate,
we can choose finite $S,S'\subseteq T$ with ${U'}\subseteq S,S'$ and $R|S$ and $R'|S'$ nondegenerate.
By Proposition \ref{22no08a}, the orbits $\GL(V)\cdot R|S$ and $\GL(V)\cdot R'|S'$ are closed.
In other words, $\GL(V)\cdot R|S=C_S$ and $\GL(V)\cdot R'|S'=C_{S'}$.
Hence, by \rf{29de14a}, $\GL(V)\cdot R|{U'}=C_{U'}=\GL(V)\cdot R'|{U'}$.
\bx

\sectz{Final remarks}

With the methods of [18] one may derive from
Theorem \ref{13me14a} the following
alternative characterization of traces of tensor representations, in terms of the exponential
rank growth of `connection matrices'
(cf.\ Freedman, Lov\'asz, and Schrijver [5]).
To this end, define, for any $f:\GG\to\oF$ and $k\in\oZ_+$, the $\WW_{k,k}\times\WW_{k,k}$ matrix $M_{f,k}$ by
$$
(M_{f,k})_{W,X}:=f(W\cdot X),
$$
for $W,X\in\WW_{k,k}$.
Then for any $f:\GG\to\oF$ and $n\in\oZ_+$:
\dy{13me14ax}{
$f$ is the trace of an $n$-dimensional tensor representation if and only if
$f(\emptyset)=1$, $f(\loop)=n$, and $\rank(M_{f,k})\leq n^{2k}$ for each $k$.
}

If $T$ is finite (which is the case in most of the examples given, and also the group
example can be described by a finite $T$ if the group is finitely generated), then all
$n$-dimensional traces form a variety, namely the closed orbit space
$\RR_T//\GL(V)$, where $V$ is $n$-dimensional.
This is a direct consequence of the fact that $p(\oF\GG)=\OO(R_T)^{\GL(V)}$ (cf.\ \rf{30de14a}).

Let us finally remark that most results of this paper have an analogue for undirected graphs, by
replacing $\GL(V)$ by the orthogonal group $O(V)$ (with respect to some nondegenerate bilinear form
on $V$).

\section*{References}\label{REF}
{\small
\begin{itemize}{}{
\setlength{\labelwidth}{4mm}
\setlength{\parsep}{0mm}
\setlength{\itemsep}{0mm}
\setlength{\leftmargin}{5mm}
\setlength{\labelsep}{1mm}
}
\item[\mbox{\rm[1]}] I. Arad, Z. Landau, 
Quantum computation and the evaluation of tensor networks,
{\em {SIAM} Journal on Computing} 39 (2010) 3089--3121. 

\item[\mbox{\rm[2]}] M. Brion, 
Introduction to actions of algebraic groups,
{\em Les cours du C.I.R.M.} 1 (2010) 1--22.

\item[\mbox{\rm[3]}] S. Chmutov, S. Duzhin, J. Mostovoy, 
{\em Introduction to Vassiliev Knot Invariants},
Cambridge University Press, Cambridge, 2012.

\item[\mbox{\rm[4]}] J. Draisma, D. Gijswijt, L. Lov\'asz, G. Regts, A. Schrijver, 
Characterizing partition functions of the vertex model,
{\em Journal of Algebra} 350 (2012) 197--206.

\item[\mbox{\rm[5]}] M.H. Freedman, L. Lov\'asz, A. Schrijver, 
Reflection positivity, rank connectivity, and homomorphisms of graphs,
{\em Journal of the American Mathematical Society} 20 (2007) 37--51.

\item[\mbox{\rm[6]}] R. Goodman, N.R. Wallach, 
{\em Symmetry, Representations, and Invariants},
Springer, Dordrecht, 2009.

\item[\mbox{\rm[7]}] P. de la Harpe, V.F.R. Jones, 
Graph invariants related to statistical mechanical models:
examples and problems,
{\em Journal of Combinatorial Theory, Series B} 57 (1993) 207--227.

\item[\mbox{\rm[8]}] T. Huckle, K. Waldherr, T. Schulte-Herbr\"uggen, 
Computations in quantum tensor networks,
{\em Linear Algebra and Its Applications} 438 (2013) 750--781. 

\item[\mbox{\rm[9]}] L.H. Kauffman, 
Knots, abstract tensors and the Yang-Baxter equation,
in: {\em Knots, Topology and Quantum Field Theories}
(Florence, 1989),
World Scientific, River Edge, N.J., 1989, pp. 179--334.

\item[\mbox{\rm[10]}] L.H. Kauffman, 
Virtual knot theory,
{\em European Journal of Combinatorics} 20 (1999) 663--690.

\item[\mbox{\rm[11]}] L.H. Kauffman, 
Introduction to virtual knot theory,
{\em Journal of Knot Theory and Its Ramifications} 21 (2012), no. 13, 1240007, 37 pp.

\item[\mbox{\rm[12]}] G. Kuperberg, 
Involutory Hopf algebras and 3-manifold invariants,
{\em International Journal of Mathematics} 2 (1991) 41--66. 

\item[\mbox{\rm[13]}] J.M. Landsberg, 
{\em Tensors: Geometry and Applications},
American Mathematical Society, Providence, R.I., 2012.

\item[\mbox{\rm[14]}] V.O. Manturov, D.P. Ilyutko, 
{\em Virtual Knots --- The State of the Art},
World Scientific, River Edge, N.J., 2013.

\item[\mbox{\rm[15]}] I.L. Markov, Y. Shi, 
Simulating quantum computation by contracting tensor networks,
{\em {SIAM} Journal on Computing} 38 (2008) 963--981. 

\item[\mbox{\rm[16]}] A. Pellionisz, R. Llin\'as, 
Tensor network theory of the metaorganization of functional geometries in the central nervous system,
{\em ￼Neuroscience} 16 (2) (1985) 245--273.

\item[\mbox{\rm[17]}] R. Penrose, 
Applications of negative dimensional tensors,
in: {\em Combinatorial Mathematics and Its Applications}
(D.J.A. Welsh, ed.),
Academic Press, London, 1971, pp. 221--244.

\item[\mbox{\rm[18]}] A. Schrijver, 
Characterizing partition functions of the vertex model by rank growth,
preprint, 2012,
ArXiv \url{http://arxiv.org/abs/1211.3561}

\item[\mbox{\rm[19]}] N. Schuch, M.M. Wolf, F. Verstraete, J.I. Cirac, 
Computational complexity of projected entangled pair states,
{\em Physical Review Letters} 98 (2007) no. 14, 140506, 4 pp. 

\item[\mbox{\rm[20]}] S. Singh, R.N.C. Pfeifer, G. Vidal, 
Tensor network decompositions in the presence of a global symmetry,
{\em Physical Review} A (3) 82 (2010) 050301, 4pp.

\item[\mbox{\rm[21]}] B. Szegedy, 
Edge coloring models and reflection positivity,
{\em Journal of the American Mathematical Society}
20 (2007) 969--988.

\item[\mbox{\rm[22]}] P. Vogel, 
Algebraic structures on modules of diagrams,
{\em Journal of Pure and Applied Algebra} 215 (2011) 1292--1339.

\end{itemize}
}

\end{document}